\documentclass[11pt]{article}
\usepackage{amsfonts}
\usepackage{epsfig}

\def\be{\begin{equation}}
\def\ee{\end{equation}}
\def\bea{\begin{eqnarray}}
\def\eea{\end{eqnarray}}
\def\bes{\begin{eqnarray*}}
\def\ees{\end{eqnarray*}}
\def\pmatrix{\left(\begin{array}{cc}}
\def\endpmatrix{\end{array}\right)}

\def\det{{\rm det}}
\def\span{{\rm span}}
\def\exp {{\rm exp}}

\def\Sp{{\rm Sp}}

\def\nn{\nonumber}
\def\<{\langle}
\def\>{\rangle}
\def\lb{\label}
\def\bs{\setminus}

\def\sp{Sp}
\def\p{\mathcal{P}}
\def\R{{\bf R}}
\def\C{{\bf C}}
\def\Z{{\bf Z}}
\def\N{{\bf N}}
\def\U{{\bf U}}
\def\Q{{\bf Q}}
\def\T{{\bf T}}

\def\ga{{\gamma}}
\def\Ga{{\Gamma}}

\def\om{{\omega}}
\def\Om{{\Omega}}

\def\lm{{\lambda}}

\def\Sg{{\Sigma}}

\def\P{{\cal P}}
\def\J{{\cal J}}

\def\I{{\cal I}}

\def\td#1{\tilde{#1}}
\def\hb{\vrule height0.18cm width0.14cm $\,$}

\def\td#1{\tilde{#1}}

\def\span{{\rm span}}
\title{ Iteration theory of $L$-index and Multiplicity of brake orbits}
\author{Chungen Liu \thanks{Partially supported by the NSF of China
973 Program of MOST. E-mail:
liucg@nankai.edu.cn}\qquad and \qquad Duanzhi Zhang
\thanks{Partially supported by NSF of China
grant 10801078.
E-mail: zhangdz@nankai.edu.cn}\\ \\
School of Mathematics and LPMC, Nankai University\\
Tianjin 300071, People's Republic of China}
\date{}
\begin{document}
\maketitle  \begin{abstract} In this paper, we first establish the
Bott-type iteration formulas  and some  abstract precise iteration
formulas of the Maslov-type index theory associated with a
Lagrangian subspace for symplectic paths. As an application, we
prove that there exist at least $\left[\frac{n}{2}\right]+1$
geometrically distinct brake orbits on every $C^2$ compact convex
symmetric hypersurface $\Sg$ in $\R^{2n}$ satisfying the reversible
condition $N\Sg=\Sg$, furthermore, if all brake orbits on this
hypersurface
 are nondegenerate, then there are at least $n$ geometrically
distinct brake orbits on it. As a consequence, we show that there
exist at least $\left[\frac{n}{2}\right]+1$ geometrically distinct
brake orbits in every bounded convex symmetric domain in $\R^{n}$,
furthermore, if all brake orbits in this domain
 are nondegenerate, then there are at least $n$ geometrically
distinct brake orbits in it. In the symmetric  case, we give a
positive answer to the Seifert conjecture of 1948 under a generic
condition.
\end {abstract}

\renewcommand{\theequation}{\thesection.\arabic{equation}}

\setcounter{equation}{0}
\section{Introduction }

Our aim of this paper is twofold. We first establish an iteration
theory of the Maslov-type index  associated with a Lagrangian
subspace of $(\mathbb{R}^{2n},\omega_0)$ for symplectic paths
starting from identity. The Bott-type iteration formulas and some
abstract precise iteration formulas are obtained here. Then as the
application of this theory, we consider the brake orbit problem on a
fixed energy hypersurface of the autonomous Hamiltonian systems. The
multiplicity results are obtained in this paper.

\subsection {Main results for the brake orbit problem}

 Let $V\in C^2(\R^n, \R)$ and $h>0$ such that $\Om\equiv \{q\in
\R^n|V(q)<h\}$ is nonempty, bounded, open and connected. Consider
the following fixed energy problem of the second order autonomous
Hamiltonian system
 \bea && \ddot{q}(t)+V'(q(t))=0, \quad {\rm for}\;q(t)\in\Om, \lb{1.1}\\
&& \frac{1}{2}|\dot{q}(t)|^2+V(q(t))= h, \qquad\forall t\in\R, \lb{1.2}\\
&& \dot{q}(0)=\dot{q}(\frac{\tau}{2})=0, \lb{1.3}\\
&& q(\frac{\tau}{2}+t)=q(\frac{\tau}{2}-t),\qquad q(t+\tau)=q(t),
    \quad \forall t\in\R.  \lb{1.4}
\eea

 A solution $(\tau,q)$ of (\ref{1.1})-(\ref{1.4}) is called a {\it
brake orbit$\;$} in $\Om$. We call two brake orbits $q_1$ and
$q_2:\R\to\R^n$ {\it geometrically distinct}  if $q_1(\R)\neq
q_2(\R)$.

We denote by $\mathcal{O}(\Om)$ and $\td{\mathcal{O}}(\Om)$ the sets
of all brake orbits and geometrically distinct brake orbits in $\Om$
respectively.

Let $J=\left(\begin{array}{cc}0&-I\\I&0\end{array}\right)$ and
$N=\left(\begin{array}{cc}-I&0\\0&I\end{array}\right)$ with $I$
being the identity in $\R^n$. Suppose that $H\in
C^2(\R^{2n}\bs\{0\},\R)\cap C^1(\R^{2n},\R)$ satisfying \be
H(Nx)=H(x),\qquad \forall\, x\in\R^{2n}.\lb{1.5}\ee

We consider the following fixed energy problem
 \bea
\dot{x}(t) &=& JH'(x(t)), \lb{1.6}\\
H(x(t)) &=& h,   \lb{1.7}\\
x(-t) &=& Nx(t),  \lb{1.8}\\
 x(\tau+t) &=& x(t),\; \forall\,t\in\R. \lb{1.9} \eea

A solution $(\tau,x)$ of (\ref{1.6})-(\ref{1.9}) is also called a
{\it brake orbit} on  $\Sg:=\{y\in\R^{2n}\,|\, H(y)=h\}$.

\noindent{\bf Remark 1.1.} It is well known that via \be H(p,q) =
{1\over 2}|p|^2 + V(q), \lb{1.10}\ee  $x=(p,q)$ and $p=\dot q$, the
elements in $\mathcal{O}(\{V<h\})$ and the solutions of
(\ref{1.6})-(\ref{1.9}) are one to one correspondent.

In more general setting, let $\Sg$ be a $C^2$ compact hypersurface
in $\R^{2n}$ bounding a compact set $C$ with nonempty interior.
Suppose $\Sg$ has non-vanishing Guassian curvature  and satisfies
the reversible condition $N(\Sg-x_0)=\Sg-x_0:=\{x-x_0|x\in \Sg\}$
for some $x_0\in C$. Without loss of generality, we may assume
$x_0=0$. We denote the set of all such hypersurface in $\R^{2n}$ by
$\mathcal{H}_b(2n)$. For $x\in \Sg$, let $N_\Sg(x)$ be the unit
outward
 normal  vector at $x\in \Sg$. Note that here by the reversible
condition there holds $N_\Sg(Nx)=NN_\Sg(x)$. We consider the
dynamics problem of finding $\tau>0$ and an absolutely continuous
curve $x:[0,\tau]\to \R^{2n}$ such that
  \bea  \dot{x}(t)&=&JN_\Sg(x(t)), \qquad x(t)\in \Sg,\lb{1.11}\\
      x(-t) &=& Nx(t), \qquad x(\tau+t) = x(t),\qquad {\rm
  for\;\; all}\;\; t\in \R.\lb{1.12}\eea

A solution $(\tau,x)$ of the problem (\ref{1.11})-(\ref{1.12}) is
a special closed characteristic on $\Sg$, here we still call it a
brake orbit on $\Sg$.

We also call two brake orbits $(\tau_1, x_1)$ and $(\tau_2,x_2)$
{\it geometrically distinct} if $x_1(\R)\ne x_2(\R)$, otherwise we
say they are equivalent. Any two equivalent brake orbits are
geometrically the same. We denote by ${\mathcal{J}}_b(\Sg)$ the set
of all brake orbits on $\Sg$, by $[(\tau,x)]$ the equivalent class
of $(\tau,x)\in {\mathcal{J}}_b(\Sg)$ in this equivalent relation
and by $\td{\mathcal{J}}_b(\Sg)$ the set of $[(\tau,x)]$ for all
$(\tau,x)\in {\mathcal{J}}_b(\Sg)$. From now on, in the notation
$[(\tau,x)]$ we always assume $x$ has minimal period $\tau$. We also
denote by $\tilde {\mathcal {J}}(\Sg)$ the set of all geometrically
distinct closed characteristics on $\Sg$.

\noindent{\bf Remark 1.2.} Similar to the closed characteristic
case, $^{\#}\td{\mathcal{J}}_b(\Sg)$ doesn't depend on the choice of
the  Hamiltonian function $H$ satisfying (\ref{1.5}) and the
conditions that $H^{-1}(\lm)=\Sg$ for some $\lm\in\R$ and $H'(x)\neq
0$ for all $x\in \Sg$.

Let $(\tau,x)$ be a solution of (\ref{1.6})-(\ref{1.9}). We consider
the boundary value problem of the linearized Hamiltonian system \bea
&&\dot{y}(t) = JH''(x(t))y(t),  \lb{1.13}\\
&&y(t+\tau)=y(t), \quad y(-t)=Ny(t), \qquad \forall t\in\R.
\lb{1.14} \eea

Denote by $\ga_x(t)$ the fundamental solution of the system
(\ref{1.13}), i.e., $\ga_x(t)$ is the solution of the following
problem
 \bea
\dot{\ga_x}(t) &=& JH''(x(t))\ga_x(t), \lb{1.15}\\
\ga_x(0) &=& I_{2n}.  \lb{1.16} \eea
 We call $\ga_x\in C([0,\tau/2],\Sp(2n))$ the {\it
associated symplectic path} of $(\tau, x)$.

The eigenvalues of $\ga_x(\tau)$ are called {\it Floquet
multipliers} of $(\tau,x)$. By Proposition I.6.13 of Ekeland's book
\cite{Ek}, the Floquet multipliers of $(\tau,x)\in
\mathcal{J}_b(\Sg)$ do not depend on the particular choice of the
Hamiltonian function $H$ satisfying conditions in Remark 1.2.

\noindent{\bf Definition 1.1.} {\it A brake orbit $(\tau,x)\in
{\mathcal{J}}_b(\Sg)$ is called  nondegenerate if 1 is its double
Floquet multiplier.}

Let $B^n_1(0)$ denote the open unit ball $\R^n$ centered at the
origin $0$. In \cite{Se1} of 1948, H. Seifert proved
$\td{\mathcal{O}}(\Om)\neq \emptyset$ provided $V'\neq 0$ on
$\partial \Om$, $V$ is analytic and $\Om$ is homeomorphic to
$B^n_1(0)$. Then he proposed his famous conjecture: {\it
$^{\#}\tilde{\mathcal{O}}(\Om)\geq n$  under the same conditions}.

After 1948, many studies have been carried out for the brake orbit
problem. S. Bolotin proved first in \cite{Bol}(also see \cite{BolZ})
of 1978 the existence of brake orbits in general setting.
 K. Hayashi in \cite {Ha1}, H. Gluck and W. Ziller
 in \cite{GZ1}, and V. Benci in \cite {Be1} in 1983-1984
  proved $^{\#}\td{\mathcal{O}}(\Om)\geq 1$ if $V$ is
$C^1$, $\bar{\Om}=\{V\leq h\}$ is compact, and $V'(q)\neq 0$ for all
$q\in \partial{\Om}$. In 1987, P. Rabinowitz in \cite{Ra1} proved
that if $H$ satisfies (\ref{1.5}), $\Sg\equiv H^{-1}(h)$ is
star-shaped, and $x\cdot H'(x)\neq 0$ for all $x\in \Sg$, then
$^{\#}\td{\mathcal{J}}_b(\Sg)\geq 1$. In 1987, V. Benci and F.
Giannoni gave a different proof of the existence of one brake orbit
in \cite{BG}.

In 1989, A. Szulkin in \cite{Sz} proved that
$^{\#}\td{\J_b}(H^{-1}(h))\geq n$, if $H$ satisfies conditions in
\cite{Ra1} of Rabinowitz and the energy hypersurface $H^{-1}(h)$ is
$\sqrt{2}$-pinched. E. van Groesen in \cite{Gro} of 1985 and A.
Ambrosetti, V. Benci, Y. Long in \cite{ABL1} of 1993 also proved
$^{\#}\td{\mathcal{O}}(\Om)\geq n$ under different pinching
conditions.

Note that the above mentioned results on the existence of multiple
brake orbits are based on certain pinching conditions. Without
pinching condition, in \cite{LZZ} Y. Long, C. Zhu and the second
author of this paper proved the following result:
 {\it For $n\ge 2$, suppose $H$ satisfies

(H1) (smoothness) $H\in C^2(\R^{2n}\bs\{0\},\R)\cap
C^1(\R^{2n},\R)$,

(H2) (reversibility) $H(Ny)=H(y)$ for all $y\in\R^{2n}$.

(H3) (convexity) $H''(y)$ is positive definite for all
$y\in\R^{2n}\bs\{0\}$,

(H4) (symmetry) $H(-y)=H(y)$ for all $y\in\R^{2n}$.

\noindent Then for any given $h>\min \{ H(y)|\; y\in \R^{2n}\}$ and
$\Sg=H^{-1}(h)$, there holds } $$^{\#}\td{\J}_b(\Sg)\ge 2.$$

As a consequence they also proved that: {\it For $n\geq 2$, suppose
$V(0)=0$, $V(q)\geq 0$, $V(-q)=V(q)$ and $V''(q)$ is positive
definite for all $q\in \R^n\bs\{0\}$. Then for $\Om\equiv
\{q\in\R^n|V(q)<h\}$ with $h>0$, there holds}
$$^{\#}\td{\mathcal{O}}(\Om)\ge 2.$$

\noindent{\bf Definition 1.2.} {\it We denote
$$\begin{array}{ll}\mathcal{H}_b^{c}(2n)=\{\Sg\in \mathcal{H}_{b}(2n)|\;\Sg\; { is\;
strictly\; convex\;} \},\\\mathcal{H}_b^{s,c}(2n)=\{\Sg\in
\mathcal{H}_{b}^c(2n)| \;-\Sg=\Sg\}.\end{array}$$}

\noindent{\bf Definition 1.3.} {\it For
$\Sg\in\mathcal{H}_b^{s,c}(2n)$, a brake orbit $(\tau,x)$ on $\Sg$
 is called symmetric if
$x(\R)=-x(\R)$. Similarly, for a $C^2$ convex symmetric bounded
domain $\Omega\subset \R^n$, a brake orbit $(\tau,q)\in \mathcal
{O}(\Omega)$ is called symmetric if $q(\R)=-q(\R)$.}

Note that a brake orbit $(\tau,x)\in \mathcal {J}_b(\Sg)$ with
minimal period $\tau$
 is symmetric if
$x(t+\tau/2)=-x(t)$ for $t\in \R$, a brake orbit $(\tau,q)\in
\mathcal {O}(\Omega)$ with minimal period $\tau$ is  symmetric if
$q(t+\tau/2)=-q(t)$ for $t\in\R$.

In this paper, we denote by $\N$, $\Z$, $\Q$ and $\R$ the sets of
positive integers, integers, rational numbers and real numbers
respectively. We denote by $\langle \cdot,\cdot\rangle$ the standard
inner product in $\R^n$ or $\R^{2n}$, by $(\cdot,\cdot)$ the inner
product of corresponding Hilbert space. For any $a\in \R$, we denote
$E(a)=\inf\{k\in \Z|k\ge a\}$ and $[a]=\sup\{k\in \Z|k\le a\}$.

The following are the main results for brake orbit problem of this
paper.

\noindent{\bf Theorem 1.1.} {\it For any
$\Sg\in\mathcal{H}_b^{s,c}(2n)$, we have }
           $$^{\#}\td{\J}_b(\Sg)\ge \left[\frac{n}{2}\right]+1 .$$

\noindent{\bf Corollary 1.1.} {\it Suppose $V(0)=0$, $V(q)\geq 0$,
$V(-q)=V(q)$ and $V''(q)$ is positive definite for all $q\in
\R^n\bs\{0\}$. Then for any given $h>0$ and $\Om\equiv
\{q\in\R^n|V(q)<h\}$, we
  have}
$$^{\#}\td{\mathcal{O}}(\Om)\ge \left[\frac{n}{2}\right]+1.$$

\noindent{\bf Theorem 1.2.} {\it For any
$\Sg\in\mathcal{H}_b^{s,c}(2n)$, suppose that all brake orbits on
$\Sg$ are nondegenerate. Then we have
         $$^{\#}\td{\J}_b(\Sg)\ge n+\mathfrak{A}({\Sg}),$$
where $2\mathfrak{A}(\Sigma)$ is the number of  geometrically
distinct asymmetric brake orbits on $\Sg$.}

As a direct consequence of Theorem 1.2, for $\Sg\in
\mathcal{H}_b^{s,c}(2n)$, if $^{\#}\td{\J}_b(\Sg)=n$ and all brake
orbits on $\Sg$ are nondegenerate, then all $[(\tau,x)]\in
\td{\mathcal {J}}_b(\Sg)$ are symmetric. Moreover, we have the
following result.

\noindent{\bf Corollary 1.2.} {\it For $\Sg\in
\mathcal{H}_b^{s,c}(2n)$, suppose $^{\#}\tilde {\mathcal
{J}}(\Sg)=n$ and all closed characteristics on $\Sg$ are
nondegenerate. Then all the $n$ closed characteristics are symmetric
brake orbits up to a suitable translation of time.}

\noindent{\bf Remark 1.3.} We note that  $^{\#}\tilde {\mathcal
{J}}(\Sg)=n$ implies $^{\#}\tilde {\mathcal {J}}_b(\Sg)\le n$, and
Theorem 1.2 implies $^{\#}\tilde {\mathcal {J}}_b(\Sg)\ge n$. So we
have $^{\#}\tilde {\mathcal {J}}_b(\Sg)=n$. Thus Corollary 1.2
follows from Theorem 1.2. Motivated by Corollary 1.2, we tend to
believe that if $\Sg\in\mathcal {H}_b^c$ and $^\#\tilde{\mathcal
{J}}(\Sg) <+\infty$, then all of them are brake orbits up to a
suitable translation of time. Furthermore, if $\Sg\in\mathcal
{H}_b^{s,c}$ and $^\#\tilde{\mathcal {J}}(\Sg)<+\infty$, then we
believe that all of them are symmetric brake orbits up to a suitable
translation of time.

\noindent{\bf Corollary 1.3.} {\it Under the same conditions of
Corollary 1.1 and the condition that all brake orbits in $\Om$ are
nondegenerate, we have
         $$^{\#}\td{\mathcal{O}}(\Om)\ge n+\mathfrak{A}(\Omega),$$
where $2\mathfrak{A}(\Omega)$ is the number of  geometrically
distinct asymmetric brake orbits in $\Om$. Moreover, if the second
order system (\ref{1.1})-(\ref{1.2}) possesses exactly $n$
geometrically distinct periodic solutions in $\Om$ and all periodic
solutions in $\Om$ are nondegenerate, then all of them are symmetric
brake orbits. }

A typical example of $\Sg\in \mathcal{H}_b^{s,c}(2n)$ is the
ellipsoid $\mathcal {E}_n(r)$ defined as follows. Let
$r=(r_1,\cdots,r_n)$ with $r_j>0$ for $1\le j\le n$. Define
$$\mathcal {E}_n(r)=\left\{x=(x_1,\cdots,x_n, y_1,\cdots,y_n)\in\R^{2n}\;
\left|\;\sum_{k=1}^n\frac{x_k^2+y_k^2}{r_k^2}=1\right.\right\}.$$
 If $r_j/r_k\notin \Q$ whenever $j\ne k$, from \cite{Ek} one can see
 that there are precisely $n$ geometrically distinct symmetric brake orbits on $\mathcal
 {E}_n(r)$ and all of them are nondegenerate.

 Since the appearance of \cite{HWZ},
Hofer, among others, has popularized in many talks the following
conjecture: {\it For $n\ge 2$, $^{\#}\tilde {\mathcal {J}}(\Sg)$ is
either $n$ or $+\infty$ for any $C^2$ compact convex hypersurface
$\Sg$ in $\R^{2n}$. } Motivated by the above conjecture and the
Seifert conjecture, we tend to believe  the following statement.

\noindent{\bf Conjecture 1.1.}  {\it For any integer $n\ge 2$, there
holds
    \bea \left\{^\#\td{\mathcal{J}}_b(\Sg)|\Sg\in
    \mathcal{H}_b^{c}(2n)\right\}=\{n, \;+\infty\}.\nn\eea}

For $\Sg\in\mathcal{H}_b^{s,c}(2n)$, Theorem 1.1 supports Conjecture
1.1 for the case $n=2$  and  Theorem 1.2 supports Conjecture 1.1 for
the nondegenerate case. However, without the symmetry assumption of
$\Sg$, the estimate $^\#\td{\mathcal{J}}_b(\Sg)\ge 2$ has not been
proved yet. It seems that there are no effective methods so far to
prove Conjecture 1.1 completely.

\subsection{Iteration formulas for Maslov-type index theory associated with a Lagrangian subspace}

We observe that the problem (\ref{1.6})-(\ref{1.9}) can be
transformed to the following problem \bea
&&\dot{x}(t) = JH'(x(t)), \nn\\
&&H(x(t))= h,   \nn\\
&&x(0)\in L_0,\;\; x(\tau/2)\in L_0,  \nn
 \eea
where $L_0=\{0\}\times \R^n\subset\R^{2n}$.

  An index theory suitable for the study of this problem was
developed in \cite{Liu2} for any Lagrangian subspace $L$. In order
to prove Theorems 1.1-1.2, we need to establish an iteration theory
for this so called $L$-index theory.

 We consider a linear Hamiltonian system
 \be\dot
x(t)=JB(t)x(t),\lb{1.17}\ee
with $B\in C([0,1],
\mathcal{L}_s(\R^{2n})$, where $ \mathcal {L}(\R^{2n})$ denotes the
set of $2n\times 2n$  real matrices  and $\mathcal{L}_s(\R^{2n})$
denotes its subset of symmetric ones. It is well known that the
fundamental solution $\gamma_B$ of (\ref{1.17}) is a symplectic path
starting from the identity $I_{2n}$ in the symplectic group
$$\Sp(2n)=\{M\in \mathcal{L}(\R^{2n})| M^TJM=J\}, $$ i.e., $\gamma_B\in \mathcal{P}(2n)$ with
$$\mathcal{P}_{\tau}(2n)=\{\gamma\in C([0,\tau],\Sp(2n))| \gamma(0)=I_{2n}\}, \;{\rm and}\;\mathcal{P}(2n)=\mathcal{P}_1(2n).$$
We denote the  nondegenerate  subset of $\mathcal{P}(2n)$ by
$$\mathcal{P}^*(2n)=\{\gamma\in \mathcal{P}(2n)| \det(\ga(1)-I_{2n})\ne 0\}.$$
 In the study of periodic solutions of Hamiltonian systems, the
Maslov-type index pair $(i(\gamma),\nu(\gamma))$ of $\gamma$  was
introduced by C. Conley and E. Zehnder in \cite{CoZ}  for
$\gamma\in\mathcal{ P}^*(2n)$ with $n\ge 2$, by Y. Long and E.
Zehnder in \cite{LZe} for  $\gamma\in\mathcal {P}^*(2)$, by Long in
\cite{Long4} and C. Viterbo in \cite{V} for $\gamma\in\mathcal
{P}(2n)$. In \cite{Long0}, Long introduced the $\omega$-index which
is an index function $(i_{\omega}(\gamma),\nu_{\omega}(\gamma))\in
\Z\times \{0,1,\cdots,2n\}$ for $\omega\in \U:=\{z\in\C|\,|z|=1\}$.

In many problems related to nonlinear Hamiltonian systems, it is
necessary to study iterations of periodic solutions. In order to
distinguish two geometrically distinct periodic solutions, one way
is to study the Maslov-type indices of the iteration paths
 of the fundamental solutions of the corresponding linearized
 Hamiltonian systems. For  $\ga\in\mathcal{P}(2n)$,  we define $\;\tilde \gamma(t)=\gamma(t-j)\gamma(1)^j$, $j\le t\le
 j+1$, $j\in\N$,  and the $k$-times iteration path of $\ga$ by $\gamma^k=\td{\gamma}|_{[0,k]}$, $\forall\, k\in \N$.
  In the paper \cite{Long0} of Long, the following result was proved
 \be i(\gamma^k)=\sum_{\omega^k=1}i_{\omega}(\gamma),
 \;\;\nu(\gamma^k)=\sum_{\omega^k=1}\nu_{\omega}(\gamma). \lb{1.18}\ee
From this result, various iteration index formulas were obtained and
were
 used to study the multiplicity and stability problems related
to the nonlinear Hamiltonian systems. We refer to the book of Long
\cite{Long1} and the references therein for these topics.

In \cite{LZZ}, Y. Long, C. Zhu and the second author of this paper
studied the multiple solutions  of the brake orbit problem on a
convex hypersurface, there they introduced indices
$(\mu_1(\gamma),\nu_1(\ga))$ and $(\mu_2(\gamma),\nu_2(\ga))$ for
 symplectic path $\gamma$. Recently, the first author of this
paper in \cite{Liu2} introduced an index theory associated with a
Lagrangian subspace for symplectic paths. For a symplectic path
$\gamma\in \mathcal{P}(2n)$, and a Lagrangian subspace $L$, by
definition the $L$-index is assigned to a pair of integers
$(i_L(\gamma), \nu_L(\gamma))\in \Z\times \{0,1,\cdots, n\}$. This
index theory is suitable for studying the Lagrangian boundary value
problems ($L$-solution, for short) related to nonlinear Hamiltonian
systems. In \cite{Liu0} the first author of this paper applied this
index theory to study the $L$-solutions of some asymptotically
linear Hamiltonian systems. The indices $\mu_1(\gamma)$ and
$\mu_2(\gamma)$ are essentially special cases of the $L$-index
$i_L(\gamma)$ for Lagrangian subspaces $L_0=\{0\}\times \R^n$ and
$L_1=\R^n\times \{0\}$ respectively up to a constant $n$.

 In order to study the
brake orbit problem, it is necessary to study the iterations of the
brake orbit. In order to do this, one way is to study the
$L_0$-index of iteration path $\gamma^k$ of the fundamental solution
$\gamma$ of the linear system (\ref{1.17}) for any $k\in \N$. In
this case, the $L_0$-iteration path $\gamma^k$ of $\gamma$ is
different from that of the general periodic case mentioned above.
Its definition is given in (\ref{4.3}) and (\ref{4.4}) below.

In 1956, Bott in \cite{Bott} established the famous iteration Morse
index formulas for closed geodesics on Riemannian manifolds. For
convex Hamiltonian systems, Ekeland  developed the similar Bott-type
iteration index formulas for Ekeland index(cf.  \cite{Ek}). In 1999,
Long in the paper \cite{Long0} established the Bott-type iteration
formulas (\ref{1.18}) for Maslov-type index. In this paper,  we
establish the following Bott-type iteration formulas for the
$L_0$-index (see Theorem 4.1 below).
\newpage
\noindent{\bf Theorem 1.3.} {\it Suppose $\gamma\in\mathcal
{P}_{\tau}(2n)$, for the iteration symplectic paths $\gamma^k$
defined in (\ref{4.3})-(\ref{uvw}) below, when $k$ is odd, there
hold \be
i_{L_0}(\gamma^{k})=i_{L_0}(\gamma^1)+\sum_{i=1}^\frac{k-1}{2}i_{\omega_{k}^{2i}}(\gamma^2),\;
\nu_{L_0}(\gamma^{k})=\nu_{L_0}(\gamma^1)+\sum_{i=1}^\frac{k-1}{2}\nu_{\omega_{k}^{2i}}(\gamma^2),
\lb{1.19}\ee
 when $k$ is even, there hold
 \be
i_{L_0}(\gamma^{k})=i_{L_0}(\gamma^1)+i^{L_0}_{\sqrt{-1}}(\gamma^1)+\sum_{i=1}^{\frac{k}{2}-1}i_{\omega_{k}^{2i}}(\gamma^2),\;
\nu_{L_0}(\gamma^{k})=\nu_{L_0}(\gamma^1)+\nu^{L_0}_{\sqrt{-1}}(\gamma^1)+\sum_{i=1}^{\frac{k}{2}-1}\nu_{\omega_{k}^{2i}}(\gamma^2),
\lb{1.20}\ee where $\omega_k=e^{\pi\sqrt{-1}/k}$ and
$(i_{\omega}(\gamma),\;\nu_{\omega}(\gamma))$ is the $\omega$ index
pair of the symplectic path $\gamma$ introduced in \cite{Long0}, and
the index pair
$(i^{L_0}_{\sqrt{-1}}(\gamma^1),\nu^{L_0}_{\sqrt{-1}}(\gamma^1))$ is
defined in Section 3.}

\noindent{\bf Remark 1.4. $\;\;$(i).} Note that the types of
iteration formulas of Ekeland and (\ref{1.18}) of Long  are the same
as that of Bott while the type of our Bott-type iteration formulas
in Theorem 1.3 is  somewhat  different from theirs. In fact, their
proofs depend on the fact that the natural decomposition of the
Sobolev space under the corresponding quadratical form is
orthogonal, but the natural decomposition in our case is no longer
orthogonal under the corresponding quadratical form.  The index pair
$(i^{L_0}_{\sqrt{-1}}(\gamma^1),\nu^{L_0}_{\sqrt{-1}}(\gamma^1))$
established in this paper is an index theory associated with two
Lagrangian subspaces.

{\bf(ii).} In \cite{LZZ}, by using $\hat{\mu}_1(x)>1$ for any brake
orbit in convex Hamiltonian systems and the dual variational method
the authors proved the existence of two geometrically distinct brake
orbits on $\Sg\in\mathcal{H}_b^{s,c}(2n)$ , where $\hat{\mu}_1(x)$
is the mean $\mu_1$-index of $x$ defined in \cite{LZZ}. Based on the
Bott-type iteration formulas in Theorem 1.3, we can deal with the
brake orbit problem more precisely to obtain the existence of more
geometrically distinct brake orbits on
$\Sg\in\mathcal{H}_b^{s,c}(2n)$.

From the Bott-type formulas in Theorem 1.3, we prove the abstract
precise iteration index formula of $i_{L_0}$ in Section 5 below.

\noindent{\bf Theorem 1.4.} {\it Let $\ga\in
\mathcal{P}_{\tau}(2n)$, $\ga^k$ is defined by
(\ref{4.3})-(\ref{uvw}) below, and $M=\ga^2(2\tau)$.
 Then for every $k\in 2\N-1$, there holds
 \bea
 i_{L_0}(\gamma^k)= i_{L_0}(\gamma^1)+\frac
{k-1}{2}(i(\gamma^2)+S^+_M(1)-C(M))
+\sum_{\theta\in(0,2\pi)}E\left(\frac{k\theta}{2\pi}\right)S_M^-(e^{\sqrt{-1}\theta})-C(M),\lb{1.21}\eea
where $C(M)$ is defined by
$$C(M)=\displaystyle\sum_{\theta\in(0,2\pi)}S^-_M(e^{\sqrt{-1}\theta})$$
and
$$S^{\pm}_M(\omega)=\lim_{\varepsilon\to 0+}i_{\omega exp(\pm
\sqrt{-1}\varepsilon)}(\gamma^2)-i_{\omega}(\gamma^2)$$ is the
splitting number of the symplectic matrix $M$ at $\omega$ for
$\omega\in \U$. (cf. \cite{Long0}, \cite{Long1}).

 For every $k\in 2\N$, there holds
\bea  i_{L_0}(\gamma^k)&=& i_{L_0}(\gamma^2)+\left(\frac
k2-1\right)\left(i(\gamma^2)+S^+_M(1)-C(M)\right)\nn\\&&-C(M)-\displaystyle\sum_{\theta\in(\pi,2\pi)}S^-_M(e^{\sqrt{-1}\theta})
+\sum_{\theta\in(0,2\pi)}E\left(\frac{k\theta}{2\pi}\right)S_M^-(e^{\sqrt{-1}\theta}).\lb{1.23}\eea}

Using the iteration formulas in Theorems 1.3-1.4, we establish the
common index jump theorem of the $i_{L_0}$-index for a finite
collection of symplectic paths starting from identity with positive
mean $i_{L_0}$-indices. In the following of this paper, we write
$(i_{L_0}(\gamma,k),\nu_{L_0}(\gamma,k))=(i_{L_0}(\gamma^k),\nu_{L_0}(\gamma^k))$
for any symplectic path $\gamma\in \mathcal {P}_{
 {\tau}}(2n)$ and $k\in \N$.

\noindent{\bf Theorem 1.5.} {\it Let $\ga_j\in \mathcal {P}_{
 {\tau_j}}(2n)$ for $j=1,\cdots,q$.  Let
 $M_j=\ga(2\tau_j)$, for $j=1,\cdots,q$. Suppose
          \be \hat{i}_{L_0}(\ga_j)>0, \quad
          j=1,\cdots,q.\lb{6.12}\ee
  Then there exist infinitely many $(R, m_1, m_2,\cdots,m_q)\in \N^{q+1}$ such that

  (i) $\nu_{L_0}(\ga_j, 2m_j\pm 1)=\nu_{L_0}(\ga_j)$,

  (ii) $i_{L_0}(\ga_j, 2m_j-1)+\nu_{L_0}(\ga_j,2m_j-1)=R-(i_{L_1}(\ga_j)+n+S_{M_j}^+(1)-\nu_{L_0}(\ga_j))$,

  (iii)$i_{L_0}(\ga_j,2m_j+1)=R+i_{L_0}(\ga_j)$.}

\subsection{Sketch of the proofs of Theorems 1.1-1.2}

For reader's convenience we briefly sketch the proofs of Theorems
1.1 and 1.2.

Fix a hypersurface $\Sg\in \mathcal{H}_b^{s,c}(2n)$ and suppose
$^\#\td{\mathcal{J}}_b(\Sg)<+\infty$, we will carry out the proof of
Theorem 1.1 in Section 7 below in the following three steps.

\noindent{\it Step 1.} Using the Clarke dual variational method, as
in \cite{LZZ}, the brake orbit problem is transformed to a fixed
energy problem of Hamiltonian systems whose Hamiltonian function is
defined by $H_\Sg(x)=j_\Sg^2(x)$ for any $x\in \R^{2n}$ in terms of
the gauge function $j_\Sg(x)$ of $\Sg$. By results in \cite{LZZ}
brake orbits in $\mathcal{J}_b(\Sg,2)$ (which is defined in Section
6 after (\ref{7.7})) correspond to critical points of
$\Phi_\Sg=\Phi|_{M_\Sg}$ where $M_\Sg$ and $\Phi$ are defined by
(\ref{7.10}) and (\ref{7.11}) in Section 6 below. Then in Section 6
we obtain the injection map $\phi: \N+K\to \mathcal
{V}_{\infty,b}(\Sg,2)\times \N$, where $K$ is a nonnegative integer
and the infinitely variationally visible subset
 $\mathcal{V}_{\infty,b}(\Sg,2)$ of $\td{\mathcal{J}}_b(\Sg,2)$ is defined in Section
 6
 such that

(i) For any $k\in \N+K$, $[(\tau,x)]\in
\mathcal{V}_{\infty,b}(\Sg,2)$ and $m\in \N$ satisfying
$\phi(k)=([(\tau \;,x)],m)$, there holds
              \be i_{L_0}(x^m)\le k-1\le i_{L_0}(x^m)+\nu_{L_0}(x^m)-1,\lb{1.25}\ee
where $x$ has minimal period $\tau$, and $x^m$ is the $m$-times
iteration of $x$ for $m\in \N$. We remind that we have written
$i_{L_0}(x)=i_{L_0}(\gamma_x)$ for a brake orbit $(\tau,x)$ with
associated symplectic path $\gamma_x$.

(ii) For any $k_j\in \N+K$, $k_1<k_2$, $(\tau_j,x_j)\in
 \mathcal{J}_b(\Sg,2)$ satisfying $\phi(k_j)=([(\tau_j \;,x_j)],m_j)$ with
$j=1,2$ and $[(\tau_1 \;,x_1)]=[(\tau_2 \;,x_2)]$, there holds
       $$m_1<m_2.$$

\noindent{\it Step 2.} Any symmetric
$(\tau,x)\in\mathcal{J}_b(\Sg,2) $ with minimal period $\tau$
satisfies
 \be
x(t+\frac{\tau}{2})=-x(t),\qquad \forall t\in \R,\lb{1.26}\ee any
asymmetric $(\tau,x)\in \mathcal{J}_b(\Sg,2)$ satisfies
 \be (i_{L_0}(x^m),\nu_{L_0}(x^m))=(i_{L_0}((-x)^m),\nu_{L_0}((-x)^m)),\quad
 \forall m\in \N.\lb{1.27}\ee
Denote the numbers of symmetric and asymmetric
 elements in $\td{\mathcal{J}}_b(\Sg,2)$ by $p$ and $2q$.
 We can
write
     $$\td{\mathcal{J}}_b(\Sg,2)=\{[(\tau_j,x_j)]|j=1,2,\cdots,p\}\cup \{[(\tau_k,x_k)],[(\tau_k,-x_k)]|k=p+1,p+2,\cdots,p+q\},$$
where $\tau_j$ is the minimal period of $x_j$ for
$j=1,2,\cdots,p+q$.

 Applying Theorem 1.5  to the associated symplectic paths of
$$(\tau_1,x_1),(\tau_2,x_2),\cdots,(\tau_{p+q},x_{p+q}),(2\tau_{p+1},x^2_{{p+1}}),(2\tau_{p+2},x^2_{{p+2}}),\cdots,(2\tau_{p+q},x^2_{{p+q}})$$
we obtain an integer $R$ large enough and the iteration times
$m_1,m_2,\cdots,m_{p+q},m_{p+q},m_{p+q+1},\cdots,m_{p+2q}$ such that
the precise information on the $(\mu_1,\nu_1)$-indices of
$(\tau_j,x_j)$'s are given in (\ref{8.49})-(\ref{8.56}).

 By the injection map $\phi$ and Step 2,
without loss of generality, we can further set
        \be \phi(R-s+1)=([(\tau_{k(s)},x_{(k(s)})],m(s))\quad {\rm
        for}\; s=1,2,\cdots,\left[\frac{n}{2}\right]+1,\lb{1.28}\ee
where $m(s)$ is the iteration time of $(\tau_{k(s)},x_{k(s)})$.

\noindent{\it Step 3.} Let \be S_1=\left\{\left.s\in
\{1,2,\cdots,\left[\frac{n}{2}\right]+1\}\right|k(s)\le
p\right\},\quad
S_2=\left\{1,2,\cdots,\left[\frac{n}{2}\right]+1\right\}\setminus{S_1}.\lb{1.29}\ee
In Section 7 we should show that
 \be ^\#S_1\le p\quad {\rm and }\quad ^\#S_2\le 2q.\lb{1.30}\ee
In fact, (\ref{1.30}) implies Theorem 1.1.

 To prove the first estimate in (\ref{1.30}),
in Section 7  below we prove the following result.

 \noindent{\bf
Lemma 1.1.} {\it Let $(\tau ,x)\in \mathcal{J}_b(\Sg,2)$ be
symmetric in the sense that $x(t+\frac{\tau}{2})=-x(t)$ for all
$t\in \R$ and $\ga$ be the associated symplectic path of $(\tau,x)$.
Set $M=\ga(\frac{\tau}{2})$. Then there is a continuous symplectic
path
      \be
      \Psi(s)=P(s) M P(s)^{-1}, \quad s\in [0,1]\lb{8.1}\ee
such that
      \be \Psi(0)=M,\qquad \Psi(1)=(-I_2)\diamond \tilde{M},\;\;\;\; \td{M}\in \Sp(2n-2),\lb{8.2}\ee
      \be \nu_1(\Psi(s))=\nu_1(M), \quad \nu_2(\Psi(s))=\nu_2(M),
      \quad  \forall \;s\in [0,1],\lb{8.3}\ee
 where
$P(s)=\left(\begin{array}{cc}\psi(s)^{-1}&0\\0&\psi(s)^T\end{array}\right)$
and $\psi$ is a continuous $n\times n$ matrix path with
$\det\psi(s)>0$ for all $s\in [0,1]$.}

 In other words, the
symplectic path $\gamma|_{[0,\tau/2]}$ is $L_j$-homotopic to a
symplectic path $\gamma^*$ with $\gamma^*(\tau/2)=(-I_2)\diamond
\td{M}$ for $j=0,1$(see Definition 2.6 below for the notion of
$L$-homotopic). This observation is essential in the proof of the
estimate
\be|(i_{L_0}(\ga)+\nu_{L_0}(\ga))-((i_{L_1}(\ga)+\nu_{L_1}(\ga))|\le
n-1\lb{abc}\ee
 in Lemma 7.1 for $\ga$ being the associated symplectic path of the symmetric
 $(\tau,x)\in \mathcal{J}_b(\Sg,2)$ in the sense that $x(t+\frac{\tau}{2})=-x(t)$ for all
$t\in \R$.
  We note that in the estimate of the
Maslov-type index $i(\gamma)$, the basic normal form theory usually
plays an important role such as in \cite{LZ}, while for the
$i_{L}$-index theory, only under the symplectic transformation of
$P(s)$ defined in Lemma 1.1, the index pairs
$(i_{L_0}(\ga),\nu_{L_0}(\ga))$ and $((i_{L_1}(\ga),\nu_{L_1}(\ga))$
are both invariant, so  the basic normal form theory can not be
applied directly.

 \noindent{\bf
Lemma 1.2.} {\it Let $(\tau ,x)\in \mathcal{J}_b(\Sg,2)$ be
symmetric in the sense that $x(t+\frac{\tau}{2})=-x(t)$ for all
$t\in \R$ and $\ga$ be the associated symplectic path of $(\tau,x)$.
 Then we have the estimate \be
i_{L_1}(\ga)+S_{\gamma(\tau)}^+(1)-\nu_{L_0}(\ga)\ge
\frac{1-n}{2}.\lb{abcd}\ee}
 {\bf Proof.} We set $\mathcal
{A}=i_{L_1}(\ga)+S_{\gamma(\tau)}^+(1)-\nu_{L_0}(\ga)$, and dually
$\mathcal {B}=i_{L_0}(\ga)+S_{\gamma(\tau)}^+(1)-\nu_{L_1}(\ga)$.
From (\ref{abc}), we have $|\mathcal{A}-\mathcal{B}|\le n-1$. It is
easy to see from Lemma 4.1 of \cite{LLZ} that
$\mathcal{A}+\mathcal{B}\ge 0$. So we have
$$\mathcal {A}\ge \frac{1-n}{2}.$$
\hfill\hb

 Combining the index estimate (\ref{abcd}) and  Lemma
7.3 below, we show that $m(s)=2m_{k(s)}$ for any $s\in S_1$. Then by
the injectivity of $\phi$ we obtain an injection map from $S_1$ to
$\{[(\tau_j,x_j)]|1\le j\le p\}$ and hence $^\#S_1\le p$.

Note that $i(\ga)=i_{\om}(\ga)$ for $\om=1$, so one can estimate
$i(\ga)+2S^+_{\gamma(\tau)}-\nu(\gamma)$ as in Lemma 4.1 of
\cite{LLZ} and $\rho_n(\Sg)$  as in \cite{LZ} by using the splitting
number theory. While the relation between the splitting number
theory and the $i_{L}$-index theory is not clear, so we have to
estimate $\mathcal{A}$ by the above method indirectly.

 To prove the second estimate of (\ref{1.30}),
using the precise index information in (\ref{8.49})-(\ref{8.56}) and
Lemmas 7.2-7.3 we can conclude that $m(s)$ is either $2m_{k(s)}$ or
$2m_{k(s)}-1$ for $s\in S_2$. Then by the injectivity of $\phi$ we
can define a map from $S_2$ to $\Ga\equiv \{[(\tau_j,x_j)]|p+1\le
j\le p+q\}$ such that any element in $\Ga$ is the image of at most
two elements in $S_2$. This yields that $^\#S_2\le 2q$.

In the following we sketch the proof of Theorem 1.2 briefly.

Suppose $^\#\td{\mathcal{J}}_b(\Sg)<+\infty$, we set \be
\td{\mathcal{J}}_b(\Sg,2)=\{[(\tau_j,x_j)]|j=1,2,\cdots,p\}\cup
\{[(\tau_k,x_k)],[(\tau_k,-x_k)]|k=p+1,p+2,\cdots,p+q\},\lb{1.31}\ee
where we have set $q=\mathfrak{A}(\Sg)$, and $\tau_j$ is the minimal
period of $x_j$ for $j=1,2,\cdots,p+q$.

Set $r=p+q$. Applying Theorem 1.5 to the associated symplectic paths
of $(\tau_1,x_1),\cdots,(\tau_r,x_r)$, we obtain an integer $R$
large enough and the iteration times $m_1,\cdots,m_r$ such that the
$i_{L_0}$-indices of iterations of $(\tau_j,x_j)$'s are given in
(\ref{9.1})-(\ref{9.3}).

Similar to (\ref{1.28}) we can set \be
\phi(R-s+1)=([(\tau_{k(s)},x_{k(s)})],m(s)) \quad{\rm for}\;
s=1,2,\cdots,n,\lb{1.32}\ee where $m(s)$ is the iteration time of
$(\tau_{k(s)},x_{k(s)})$. Then by Lemma 7.3,
(\ref{9.1})-(\ref{9.3}), and that $x_j^m$ is nondegenerate for $1\le
j\le r$ and $m\in\N$ , we prove that $m(s)=2m_{k(s)}$. Then by the
injectivity of $\phi$ we have
$$^\#\td{\mathcal{J}}_{b}(\Sg)=^\#\td{\mathcal{J}}_{b}(\Sg,2)=p+2q=r+q\ge
n+q=n+\mathfrak{A}(\Sg).$$

 This paper is organized as follows. In
Section 2, we briefly introduce the $L$-index theory associated with
Lagrangian subspace $L$ for symplectic paths and give upper bound
estimates for $|i_{L_0}-i_{L_1}|$ and
$|(i_{L_0}+\nu_{L_0})-(i_{L_1}+\nu_{L_1})|$. In Section 3, we
introduce an $\om$-index theory for symplectic paths associated with
a Lagrangian subspace. Then in Section 4 we establish the Bott-type
iteration formulas of the Maslov-type indices $i_{L_0}$ and
$i_{L_1}$. Based on these Bott-type iteration formulas we prove
Theorems 1.4 and 1.5 in Section 5.  In Section 6, we obtain the
injection map $\phi$ which is also  basic in the proofs of Theorems
1.1 and 1.2. Based on these results in Sections 5 and 6, we prove
Theorem 1.1 in Section 7, and we finally prove Theorem 1.2 in
Section 8.

\setcounter{equation}{0}
\section {Maslov type $L$-index theory associated with a Lagrangian subspace for symplectic paths
} 

 In this section, we give a brief introduction to the  Maslov type $L$-index theory.
  We refer to the papers \cite{Liu2}  and \cite{Liu0} for the
 details.

  Let $(\R^{2n}, \omega_0)$ be the standard linear symplectic space with $\omega_0=\sum_{j=1}^n dx_j\wedge dy_j$.
 A Lagrangian subspace $L$ of $(\R^{2n}, \omega_0)$ is an $n$ dimensional subspace
 satisfying $\omega_0|_L=0$. The set of all Lagrangian subspaces in $(\R^{2n}, \omega_0)$ is denoted by $\Lambda(n)$.

 For a symplectic path $\gamma\in \mathcal{P}(2n)$, we write it in the following
 form
\be\gamma(t)=\left(\begin{array}{cc}S(t)&V(t)\\T(t)&U(t)\end{array}\right),\lb{2.1}\ee
 where $S(t), T(t), V(t), U(t)$ are $n\times n$ matrices.
 The $n$ vectors coming from the columns of the matrix $\left(\begin{array}{c}
 V(t)\\U(t)\end{array}\right)$ are linear independent and they span a
 Lagrangian subspace path of $(\R^{2n}, \omega_0)$.
 For  ${L_{0}}=\{0\}\times \R^n\in \Lambda(n)$,
 we define the following two subsets of $\Sp(2n)$ by
 $$\Sp(2n)_{L_{0}}^*=\{M\in\Sp(2n)|\, \det V\neq 0\},$$
 $$\Sp(2n)_{L_{0}}^0=\{M\in\Sp(2n)|\, \det V= 0\},$$
 for  $M=\left(\begin{array}{cc}S & V\\T & U\end{array}\right)$.

 Since the space $\Sp(2n)$ is path connected, and the set of $n\times n$
 non-degenerate matrices  has two path connected components consisting of matrices with positive and
  negative determinants respectively.
 We denote by
 $$\Sp(2n)_{L_{0}}^{\pm}=\{M\in\Sp(2n)|\, \pm\det V>0 \}, $$
 $$\mathcal{P}(2n)_{L_{0}}^*=\{\gamma\in \mathcal{P}(2n)|\, \gamma(1)\in\Sp(2n)_{L_{0}}^*\},$$
 $$ \mathcal{P}(2n)_{L_{0}}^0=\{\gamma\in \mathcal{P}(2n)|\, \gamma(1)\in\Sp(2n)_{L_{0}}^0\}. $$

 \noindent{\bf Definition 2.1.}(\cite{Liu2}) {\it We define the ${L_{0}}$-nullity of any
 symplectic path $\gamma\in \mathcal{P}(2n)$ by
 \be\nu_{L_{0}}(\gamma)=\dim\ker V(1)
 \lb{2.2}\ee
 with the $n\times n$ matrix function $V(t)$  defined in
 (\ref{2.1}).}

 We note that the complex matrix
 $U(t)\pm\sqrt{-1}V(t)$ is invertible. We define a complex matrix function by
 \be\mathcal{Q}(t)=[U(t)-\sqrt{-1}V(t)] [U(t)+\sqrt{-1}V(t)]^{-1}. \lb{2.3}\ee
 The matrix $\mathcal {Q}(t)$ is unitary for
 any $t\in [0,1]$. We denote by
 $$M_+= \pmatrix 0 & I_n\\ -I_n & 0\endpmatrix, \;\; M_-=\pmatrix 0 & J_n\\
  -J_n & 0\endpmatrix, \;\;J_n={\rm diag}(-1,1,\cdots,1). $$
  It is clear that $M_{\pm}\in\Sp(2n)_{L_{0}}^{\pm}$.

  For a path $\gamma\in \p(2n)_{L_{0}}^*$,  we define a symplectic path by
 \be\tilde {\gamma}(t)=\left\{\begin{array}{lr} I\cos \frac{(1-2t)\pi}{2}+J\sin\frac{(1-2t)\pi}{2}, \;\;& t\in [0,1/2],\\
  \gamma(2t-1),\; & t\in [1/2,1]\end{array}\right. \lb{2.4}\ee
  and choose a symplectic path $\beta(t)$ in $\Sp(2n)_{L_{0}}^*$ starting from
  $\gamma(1)$ and ending at $M_+$ or $M_-$ according to $\gamma(1)\in\sp(2n)_{L_{0}}^+$
  or $\gamma(1)\in\sp(2n)_{L_{0}}^-$, respectively. We now define a joint path by
  \be\bar{\gamma}(t)=\beta*\tilde {\gamma}:=\left\{\begin{array}{lr} \tilde {\gamma}(2t), \;\;& t\in [0,1/2],\\
  \beta(2t-1),\;\; & t\in [1/2,1].\end{array}\right. \lb{2.5}\ee
  By the definition, we see that the symplectic path $\bar{\gamma}$
  starts from $-M_+$ and ends at either $M_+$ or $M_-$.
  As above, we define
  \be\bar {\mathcal{Q}}(t)=[\bar {U}(t)-\sqrt{-1}\bar {V}(t)] [\bar {U}(t)+\sqrt{-1}\bar {V}(t)]^{-1}.
  \lb{2.6}\ee
  for $\bar {\gamma}(t)=\pmatrix \bar {S}(t) & \bar {V}(t)\\\bar {T}(t) & \bar
  {U}(t)\endpmatrix$. We can choose a continuous function $\bar
  {\Delta}(t)$ on $[0,1]$ such that
  \be\det \bar {\mathcal {Q}}(t)=e^{2\sqrt{-1}\bar{\Delta}(t)}. \lb{ 2.7}\ee
  By the above arguments, we see that the number $\frac{1}{\pi}(\bar
  {\Delta}(1)-\bar{\Delta}(0))\in \Z$ and it does not depend on
  the choice of the function $\bar{\Delta}(t)$.

\noindent{\bf Definition 2.2.}(\cite{Liu2}) {\it For a symplectic
path $\gamma\in \p(2n)_{L_{0}}^*$, we define the ${L_{0}}$-index of
$\gamma$ by \be i_{L_{0}}(\gamma)=\frac{1}{\pi}(\bar
  {\Delta}(1)-\bar{\Delta}(0)). \lb {2.8}\ee}
 \noindent{\bf Definition 2.3.}(\cite{Liu2}) {\it For a symplectic path $\gamma\in
\p(2n)_{L_{0}}^0$, we define the ${L_{0}}$-index of $\gamma$ by \be
i_{L_{0}}(\gamma)=\inf\{i_{L_{0}}(\gamma^*)|\,\gamma^*\in
\p(2n)_{L_{0}}^*, { \,\gamma^*\,is\,sufficiently\,close \,to}
\,\gamma \}. \lb {2.9}\ee}
 $\quad$ In the general situation, let
$L\in \Lambda(n)$.  It is well known that $\Lambda(n)=U(n)/O(n)$,
this means that for any linear subspace $L\in \Lambda(n)$,
there is an orthogonal symplectic matrix $P=\pmatrix A & -B\\
B & A\endpmatrix$ with $A\pm \sqrt {-1}B\in U(n)$ such that
$PL_0=L$.  We define the conjugated symplectic path $\gamma_c\in
\p(2n)$ of $\gamma$ by $\gamma_c(t)=P^{-1}\gamma(t)P$.

\noindent{\bf Definition 2.4.}(\cite{Liu2}) {\it We define the
$L$-nullity of any
 symplectic path $\gamma\in \p(2n)$ by
\be\nu_L(\gamma)=\dim\ker V_c(1), \lb{2.10}\ee
  the $n\times n$ matrix function $V_c(t)$  is defined in (\ref{2.1}) with the symplectic path
 $\gamma$ replaced by $\gamma_c$, i.e.,
 \be\gamma_c(t)=\pmatrix S_c(t) & V_c(t)\\T_c(t) & U_c(t)\endpmatrix. \lb {2.11}\ee
 }

 \noindent{\bf Definition 2.5.}(\cite{Liu2}) {\it For a symplectic path $\gamma\in
\p(2n)$, we define the ${L}$-index of $\gamma$ by \be
i_{L}(\gamma)=i_{L_0}(\gamma_c). \lb {2.12}\ee } We define a Hilbert
space $E^1=E^1_{L_0}=W^{1/2,2}_{L_0}([0,1],\R^{2n})$ with $L_0$
boundary conditions by \bea E^1_{L_0}=\left\{x\in
L^2([0,1],\R^{2n})|
x(t)=\sum_{j\in\Z}\exp(j\pi tJ) \left(\begin{array}{c} 0\\
a_j\end{array}\right), a_j\in \R^n,
\;\|x\|^2:=\sum_{j\in\Z}(1+|j|)|a_j|^2<\infty\right\}.\nn\eea

For any Lagrangian subspace $L\in \Lambda(n)$, suppose $P\in
\Sp(2n)\cap O(2n)$ such that $L=PL_0$. Then we define
$E^1_L=PE^1_{L_0}$. We define two operators on $E^1_L$ by \be
(Ax,y)=\int^1_0\<-J\dot x,y\>\,dt,\;\; (Bx,y)=\int^1_0\langle B(t)
x,y\>\,dt,\;\;\forall\; x,\,y\in E^1_{L},\lb{2.13}\ee where
$(\cdot,\cdot)$ is the inner product in $E^1_L$ induced from
$E^1_{L_0}$.

By the Floquet theory we have
$$\nu_{L}(\gamma_B)=\dim\ker(A-B).  $$

We denote by $E^{L_0}_m=\left\{z\in E^1_{L_0}\left|\,
z(t)=\displaystyle\sum_{k=-m}^m-J\exp(k\pi tJ)a_k\right.\right\}$
the finite dimensional truncation of $E^1_{L_0}$, and
$E^L_m=PE^{L_0}_m$.

Let $P_m:\,E^1_L\to E^L_m$ be the orthogonal projection for
$m\in\N$. Then $\Gamma=\{P_m|\;m\in\N\}$ is a Galerkin approximation
scheme with respect to $A$ defined  in (\ref{2.13}), i.e., there
hold
$$P_m\to I\;{\rm strongly\;as}\;m\to \infty $$
and
$$P_mA=AP_m. $$

For $d>0$, we denote by $m^*_d(\cdot)$ for $*=+, 0, -$ the dimension
of the total eigenspace corresponding to the eigenvalues $\lambda$
belonging to $[d,+\infty), (-d,d)$ and $(-\infty, -d]$ respectively,
and denote by $m^*(\cdot)$  for $*=+,0,-$ the  dimension of the
total eigenspace corresponding to the eigenvalues $\lambda$
belonging to $(0,+\infty), \{0\}$ and $(-\infty, 0)$ respectively.
For any self-adjoint operator $T$, we denote $T^{\sharp}=(T|_{Im
T})^{-1}$ and $P_mTP_m=(P_mTP_m)|_{E^L_m}$.

If $\gamma_B\in \mathcal {P}(2n)$ is the fundamental solution of the
system (\ref{1.17}), we write $i_L(B)=i_L(\gamma_B)$ and
$\nu_L(B)=\nu_L(\gamma_B)$.
 The following Galerkin approximation result will be used in this
 paper.

\noindent{\bf Proposition 2.1.} (Theorem 2.1 of \cite{Liu0}) {\it
For any $B\in C([0,1], \mathcal{L}_s(\R^{2n}))$ with the $L$-index
pair $(i_L(B),\nu_L(B))$ and any constant $0<d\le \frac
14\|(A-B)^{\sharp}\|^{-1}$, there exists $m_0>0$ such that for $m\ge
m_0$, we have
\bea && m^+_d(P_m(A-B)P_m)=mn-i_L(B)-\nu_L(B),\nn\\
 && m^-_d(P_m(A-B)P_m)=mn+i_L(B)+n,\lb{2.14}\\&&
 m^0_d(P_m(A-B)P_m)=\nu_L(B).\nn\eea
}
 $\quad$ The Galerkin approximation formula for the Maslov-type
index theory associated with periodic boundary value was proved in
\cite{FQ} by Fei and Qiu.

 \noindent{\bf Remark 2.1.} Note that $mn=m^-_d(P_mAP_m)$, so we
have
 $m^-_d(P_m(A-B)P_m)-mn=I(A,A-B)$, where
 $I(A,A-B)$ is defined in
 Definition 3.1 below. So we have
 \be I(A,A-B)=i_L(B)+n.\lb{c1}\ee

\noindent{\bf Definition 2.6.} (\cite{Liu2}) {\it For two paths
$\gamma_0,\;\gamma_1\in
 \mathcal{P}(2n)$, we say that they are $L$-homotopic and denoted by
 $\gamma_0\sim_L\gamma_1$, if there is a map $\delta:[0,1]\to
 \mathcal{P}(2n)$ such that $\delta(j)=\gamma_j$ for $j=0,1$, and
 $\nu_L(\delta(s))$ is constant for $s\in [0,1]$.
 }

 For any  two $2k_i\times 2k_i$ matrices
of square block form, $M_i=\pmatrix A_i & B_i\\
                                            C_i & D_i \endpmatrix$ with $i=1, 2$,
                                             the $\diamond$-product of
$M_1$ and $M_2$ is defined  to be the $2(k_1+k_2)\times 2(k_1+k_2)$
matrix
$$ M_1\diamond M_2=\left(\begin{array}{cccc} A_1 & 0 & B_1 & 0\\
                                               0 & A_2 & 0 & B_2\\
                                           C_1 & 0 & D_1 & 0\\
                                                0 & C_2 & 0 & D_2 \end{array}\right).  $$

  \noindent{\bf Theorem 2.1.}(\cite{Liu2}) {\it If  $\gamma_0\sim_L\gamma_1$, there hold
  $$i_L(\gamma_0)=i_L(\gamma_1),\;\nu_L(\gamma_0)=\nu_L(\gamma_1).$$
}
 \noindent{\bf Theorem 2.2.}(\cite{Liu2}) {\it If $\gamma=\gamma_1\diamond \gamma_2\in
  \mathcal{P}(2n)$, and correspondingly $L=L'\oplus L''$, then
  $$i_L(\gamma)=i_{L'}(\gamma_1)+i_{L''}(\gamma_2),\;\nu_L(\gamma)=\nu_{L'}(\gamma_1)+\nu_{L''}(\gamma_2).$$
}
\noindent{\bf Theorem 2.3.} {\it For $L_0=\{0\}\times \R^n,
L_1=\R^n\times \{0\}$, then for $\gamma\in \mathcal{P}(2n)$ \be
|i_{L_0}(\gamma)-i_{L_1}(\gamma)|\le
n,\;|i_{L_0}(\gamma)+\nu_{L_0}(\gamma)-i_{L_1}(\gamma)-\nu_{L_1}(\gamma)|\le
n.\lb{2.14''}\ee Moreover,  the left hand sides of the above two
inequalities depend only on the end matrix $\gamma(1)$, in
particular, if $\gamma(1)\in O(2n)\cap Sp(2n)$, there holds \be
i_{L_0}(\gamma)=i_{L_1}(\gamma).\lb{2.14'}\ee }
 \noindent {\bf
Proof.} We only need to prove the first inequality in
 (\ref{2.14''})
\be|i_{L_0}(\gamma)-i_{L_1}(\gamma)|\le n.\lb{2.15}\ee
 For the
second inequality in (\ref{2.14''}), we can choose a symplectic path
$\gamma_1$ such that
$$i_{L_0}(\gamma)+\nu_{L_0}(\gamma)=i_{L_0}(\gamma_1),\;i_{L_1}(\gamma)+\nu_{L_1}(\gamma)=i_{L_1}(\gamma_1).$$
Then by (\ref{2.15}) we have
       \bea |i_{L_0}(\gamma_1)-i_{L_1}(\gamma_1)|\le n\nn\eea
which yields the second inequality of (\ref{2.14''}).

Note that (\ref{2.15}) holds from Theorem 3.3 of \cite{LZZ} and
Proposition 5.1 below. Here we give another proof directly from the
definitions of $i_{L_0}$ and $i_{L_1}$.

We write $\bar\gamma(t)$ in (\ref{2.5}) in its polar decomposition
form $\bar\gamma(t)=\bar O(t)\bar P(t)$, $\bar O(t)\in O(2n)\cap
Sp(2n)$, and $\bar P(t)$ is a positive definite matrix function. By
 (4.1) of \cite{Liu2} we have
$$\bar \Delta(t)=\bar \Delta_{\bar O}(t)+\bar \Delta_{\bar P}(t).$$
Since $\bar P(0)=\bar P(1)=I_{2n}$ and the set of positive definite
symplectic matrices is contractible, we have
$$\bar \Delta_{\bar P}(1)-\bar \Delta_{\bar P}(0)=0,$$
so $$\bar \Delta(1)-\bar \Delta(0)=\bar \Delta_{\bar O}(1)-\bar
\Delta_{\bar O}(0).$$ On the other hand,
$\gamma_c(t)=J^{-1}\gamma(t)J=O(t)(J^{-1}P(t)J)$. We also write
$\bar \gamma_c=\bar O_c \bar P_c$. So by the definitions of $\bar
\ga_c$ and $\bar \ga$ we have $\bar O_c(t)=\bar O(t)$ for $t\in
[0,\frac{1}{2}]$ in (\ref{2.5}). Then (\ref{2.15}) follows from the
fact that the only difference between $\bar O_c$ and $\bar O$ is
that $\td {\gamma}_c(1)$ and $\td {\gamma}(1)$ in (\ref{2.4}) may be
connected to different matrices $M^+$ or $M^-$ by $\beta_c$ and
$\beta$ in (\ref{2.5}) respectively. The statement that  the left
hand sides of the two inequalities in (\ref{2.14''}) depend only on
the end matrix $\gamma(1)$ is a consequence of Corollary 4.1 of
\cite{Liu2}. For the proof of (\ref{2.14'}), suppose $\gamma(1)\in
O(2n)\cap Sp(2n)$, we can take $\gamma(t)\in O(2n)\cap Sp(2n)$ since
the number on the left side of inequality (\ref{2.15}) depends only
on $\gamma(1)$. For $\gamma(t)\in O(2n)\cap Sp(2n)$, we have
$\gamma_c(t)=J^{-1}\gamma(t)J=\gamma(t)$. Thus we have
$i_{L_0}(\gamma)=i_{L_1}(\gamma)$. \hfill\hb

\noindent{\bf Theorem 2.4.} (Lemma 5.1 of \cite{Liu2}) {\it If
$\gamma\in \mathcal{P}(2n)$ is the fundamental solution of
$$\dot x(t)=JB(t)x(t)$$ with symmetric matrix function
$B(t)=\pmatrix b_{11}(t) & b_{12}(t)\\b_{21}(t) & b_{22}(t)
\endpmatrix$ satisfying $b_{22}(t)>0$ for any $t\in R$, then there holds
$$i_{L_0}(\gamma)=\sum_{0<s<1}\nu_{L_0}(\gamma_s),\;\gamma_s(t)=\gamma(st).$$
Similarly, if $b_{11}(t)>0$ for any $t\in \R$, there holds
$$i_{L_1}(\gamma)=\sum_{0<s<1}\nu_{L_1}(\gamma_s),\;\gamma_s(t)=\gamma(st).$$  }

\setcounter{equation}{0}
\section{ $\omega$-index theory associated with
a Lagrangian subspace for symplectic paths }

Let $E$ be a separable Hilbert space, and $Q=A-B: E\to E$ be a
bounded self-adjoint linear operators with $B:E\to E$ being a
compact self-adjoint operator. Suppose that $N=\ker Q$ and $\dim
N<+\infty$. $Q|_{N^{\bot}}$ is invertible. $P:E\to N$ is the
orthogonal projection. We denote  $d=\frac 14
\|(Q|_{N^{\bot}})^{-1}\|^{-1}$. Suppose
$\Gamma=\{P_k|k=1,2,\cdots\}$ is the Galerkin approximation sequence
of $A$ with

\quad(1)  $E_k:=P_kE$ is finite dimensional for all $k\in\N$,

\quad(2)  $P_k\to I$ strongly as $k\to +\infty$

\quad(3)  $P_kA=AP_k$.

For a self-adjoint operator $T$, we denote by $M^{*}(T)$ the
eigenspaces of $T$ with eigenvalues belonging to $(0,+\infty)$,
$\{0\}$ and $(-\infty, 0)$ with  $*=+,0$ and $*=-$, respectively. We
denote by $m^*(T)=\dim M^*(T)$. Similarly, we denote by $M_d^{*}(T)$
the $d$-eigenspaces of $T$ with eigenvalues belonging to
$(d,+\infty)$, $(-d,d)$ and $(-\infty, -d)$ with  $*=+,0$ and $*=-$,
respectively. We denote by $m_d^*(T)=\dim M_d^*(T)$.

\noindent{\bf Lemma 3.1.} {\it There exists $m_0\in \N$ such that
for all $m\ge m_0$, there hold}

\be m^-(P_m(Q+P)P_m)=m_d^-(P_m(Q+P)P_m) \lb {3.1}\ee  and \be
m^-(P_m(Q+P)P_m)=m_d^-(P_mQP_m). \lb{3.2}\ee

 \noindent{\bf Proof.} The proof of
(\ref{3.1}) is essential the same as that of Theorem 2.1 of
\cite{Fei}, we note that $\dim\ker(Q+P)=0$.

By considering the operators $Q+sP$ and $Q-sP$ for small $s>0$, for
example $s<\min \{1, d/2\}$, there exists $m_1\in\N$ such that
 \be
m^-_{d}(P_mQP_m)\le m^-(P_m(Q+sP)P_m), \;\forall \, m\ge m_1 \lb
{3.3}\ee and \be m^-_{d}(P_mQP_m)\ge
m^-(P_m(Q-sP)P_m)-m^0_{d}(P_mQP_m),\; \forall \, m\ge m_1.
\lb{3.4}\ee
 In fact, the claim (\ref{3.3}) follows from
$$P_m(Q+sP)P_m=P_mQP_m+sP_mPP_m $$
and for $x\in M^-_{d}(P_mQP_m)$,
$$(P_m(Q+sP)P_mx,x)\le -d\|x\|^2+s\|x\|^2\le -\frac{d}{2}\|x\|^2. $$
The claim (\ref{3.4}) follows from that for $x\in
M^-(P_m(Q-sP)P_m)$,
$$(P_mQP_mx,x)\le s(P_mPP_mx,x)< d\|x\|^2. $$
 By the Floquet theory, for $m\ge m_1$ we have
$m^0_{d}(P_mQP_m)=\dim N=\dim Im(P_mPP_m)$, and by
$Im(P_mPP_m)\subseteq M^0_{d}(P_mQP_m)$ we have $Im(P_mPP_m)=
M^0_{d}(P_mQP_m)$. It is easy to see that $M^0_{d}(P_mQP_m)\subseteq
M^+_{d}(P_m(Q+sP)P_m)$. By using
$$P_m(Q-sP)P_m=P_m(Q+sP)P_m-2sP_mPP_m $$
we have \be m^-(P_m(Q-sP)P_m)\ge m^-(P_m(Q+sP)P_m)+m^0_{d}(P_mQP_m),
\;\forall\,m\ge m_1.  \lb {3.5}\ee Now (\ref{3.2}) follows from
(\ref{3.3})-(\ref{3.5}). \hfill\hb

Since $M^-(Q+P)=M^-(Q)$ and the two operators $Q+P$ and $Q$ have the
same negative spectrum, moreover, $P_m(Q+P)P_m\to Q+P$ and
$P_mQP_m\to Q$ strongly, one can prove (\ref{3.2}) by the spectrum
decomposition theory.

 The following result was proved in \cite{CLL}.

\noindent{\bf Lemma 3.2.} {\it Let $B$ be a linear symmetric
compact
 operator, $P:E\to \ker A$ be the orthogonal projection. Suppose that
 $A-B$ has a bounded inverse. Then the difference of the Morse
 indices
 $$m^-(P_m(A-B)P_m)-m^-(P_m(A+P)P_m) $$
 eventually becomes a constant independent of $m$, where $A:E\to E$ is a
 bounded self-adjoint operator with a finite dimensional kernel,
 and the restriction $A|_{(\ker A)^\bot}$ is invertible, and
 $\Gamma=\{P_k\}$ is a Galerkin approximation sequence with respect to $A$.
 }

 By Lemmas 3.1 and 3.2, we have the following result.

 \noindent{\bf Lemma 3.3.} {\it Let $B$ be a linear symmetric  compact
 operator. Then the difference of the $d$-Morse
 indices
 \be m_d^-(P_m(A-B)P_m)-m_d^-(P_mAP_m) \lb {3.6}\ee
 eventually becomes a constant independent of $m$, where $d>0$ is determined by the operators $A$ and $A-B$.
 Moreover $m^0_d(P_m(A-B)P_m)$ eventually becomes a constant independent of
 $m$ and for large $m$, there holds
\be m_d^0(P_m(A-B)P_m)=m^0(A-B). \lb{3.7}\ee }
\noindent{\bf Proof.}
We only need to prove (\ref{3.7}). It is easy to show that there is
a constant $m_1>0$ such that for $m\ge m_1$
$$\dim P_m\ker (A-B)=\dim\ker(A-B). $$
Since $B$ is compact, there is $m_2\ge m_1$ such that for $m\ge m_2$
$$\|(I-P_m)B\|\le 2d. $$
Take $m\ge m_2$, let $E_m=P_m\ker(A-B)\bigoplus Y_m$, then
$Y_m\subseteq {\rm Im} (A-B)$. For $y\in Y_m$ we have
$$y=(A-B)^{\sharp} (A-B)y=(A-B)^{\sharp}(P_m(A-B)P_my+(P_m-I)By).$$
It implies
$$\|P_m(A-B)P_my\|\ge 2d\|y\|, \;\forall y\in Y_m. $$
Thus we have \be  m_d^0(P_m(A-B)P_m)\le m^0(A-B).\lb{3.8}\ee
 On the other hand,
for $x\in P_m\ker(A-B)$, there exists $y\in \ker(A-B)$, such that
$x=P_my$. Since $P_m\to I$ strongly, there exists $m_3\ge m_2$ such
that for $m\ge m_3$
$$\|I-P_m\|<\frac 12,\;\;P_m(A-B)(I-P_m)\le \frac d2. $$
So we have
$$\|P_m(A-B)P_m x\|=\|P_m(A-B)(I-P_m) y\|\le \frac d2\|y\|<d\|x\|. $$
It implies that \be m_d^0(P_m(A-B)P_m)\ge m^0(A-B).\lb{3.9}\ee
(\ref{3.7}) holds from (\ref{3.8}) and (\ref{3.9}). \hfill\hb

\noindent{\bf Definition 3.1.} {\it For the self-adjoint Fredholm
operator $A$ with a Galerkin approximation sequence $\Gamma$ and the
self-adjoint compact operator $B$ on Hilbert space $E$, we define
the relative index by \be I(A,A-B)=
m_d^-(P_m(A-B)P_m)-m_d^-(P_mAP_m),\;\;\;\; m\ge m^*, \lb {3.10}\ee
 where
$m^*>0$ is a  constant large enough such that the difference in
(\ref{3.6}) becomes a constant independent of $m\ge m^*$. }

The spectral flow for a parameter family of linear self-adjoint
Fredholm operators was introduced by Atiyah, Patodi and Singer in
\cite{APS}. The following result shows that the relative index in
Definition 3.1 is a spectral flow.

 \noindent{\bf Lemma 3.4.} {\it For  the
operators $A$ and $B$ in Definition 3.1, there holds
 \be I(A,A-B)=-{\rm sf}\{A-sB,\,0\le s\le 1\}, \lb
{3.11}\ee
 where ${\rm sf}(A-sB,\,0\le s\le 1)$ is the spectral flow of the operator family $A-sB$, $s\in[0,1]$ (cf.
\cite{ZL}).}

\noindent{\bf Proof.}  For simplicity, we set  $I_{\rm
sf}(A,A-B)=-{\rm sf}\{A-sB,\,0\le s\le 1\}$ which is exact the
relative Morse index defined in \cite{ZL}. By the Galerkin
approximation formula in Theorem 3.1 of \cite{ZL},
   \bea I_{\rm
sf}(A,A-B)=I_{\rm sf}(P_mAP_m,\,P_m(A-B)P_m)\lb{a1}\eea if
$\ker(A)=\ker(A-B)=0$.

By (2.17) of \cite{ZL}, we have \bea I_{\rm
sf}(P_mAP_m,\,P_m(A-B)P_m)&=&m^-(P_m(A-B)P_m)-m^-(P_mAP_m)\nn\\
&=&m_d^-(P_m(A-B)P_m)-m_d^-(P_mAP_m)\nn\\&=&I(A,A-B)\lb{a3}\eea
 for
$d>0$ small enough. Hence (\ref{3.11}) holds in the nondegenerate
case. In general, if $\ker(A)\ne 0$ or $\ker(A-B)\ne 0$, we can
choose $d>0$ small enough such that $\ker(A+d {\rm Id})=\ker(A-B+d
{\rm Id})=0$, here ${\rm Id}:\; E\to E$ is the identity operator. By
(2.14) of \cite{ZL} we have
 \bea I_{\rm
sf}(A,A-B)&=&I_{\rm sf}(A,A+d {\rm Id})+ I_{\rm sf}(A+d {\rm
Id},A-B+d {\rm Id})+I_{\rm sf}(A-B+d
{\rm Id},A-B)\nn\\
&=&I_{\rm sf}(A+d {\rm Id},A-B+d {\rm Id})=I(A+d {\rm Id},A-B+d {\rm Id})\nn\\
&=&m^-(P_m(A-B+d {\rm Id})P_m)-m^-(P_m(A+d
{\rm Id})P_m)\nn\\
&=&m_d^-(P_m(A-B)P_m)-m_d^-(P_mAP_m)= I(A,A-B).\lb{a2}
 \eea
In the second equality of (\ref{a2}) we note that $I_{\rm sf}(A,A+d
{\rm Id})=I_{\rm sf}(A-B+d {\rm Id},A-B)=0$ for $d>0$ small enough
since the spectrum of $A$ is discrete and $B$ is a compact operator,
in the third and the forth equalities of (\ref{a2}) we have applied
(\ref{a3}). \hfill\hb

 A similar way to define the relative index of two
operators was appeared in \cite{CLL}. A different way to study the
relative index theory was appeared in \cite{Fei}.

For $\omega=e^{\sqrt{-1}\theta}$ with $\theta\in\R$, we define a
Hilbert space $E^{\omega}=E^{\omega}_{L_0}$ consisting of those
$x(t)$ in $L^2([0,1], \C^{2n})$ such that $e^{-\theta t J}x(t)$ has
Fourier expending
$$e^{-\theta t J}x(t)=\sum_{j\in \Z}e^{j\pi tJ}\pmatrix 0\\a_j\endpmatrix,\;a_j\in \C^n  $$
with
$$\|x\|^2:=\sum_{j\in\Z}(1+|j|)|a_j|^2<\infty. $$
For $x\in E^{\omega}$, we can write \bea x(t)&=&e^{\theta
tJ}\sum_{j\in\Z}e^{j\pi tJ}\pmatrix 0\\a_j\endpmatrix
 =\sum_{j\in\Z}e^{(\theta+j\pi)tJ}\pmatrix 0\\a_j\endpmatrix \nn\\
 &=&\sum_{j\in\Z}e^{(\theta+j\pi)t\sqrt{-1}}\pmatrix \sqrt{-1}a_j/2\\a_j/2\endpmatrix+
 e^{-(\theta+j\pi)t\sqrt{-1}}\pmatrix
 -\sqrt{-1}a_j/2\\a_j/2\endpmatrix.\lb{4.7}\eea
So we can write \be x(t)=\xi(t)+N\xi(-t), \;
\xi(t)=\sum_{j\in\Z}e^{(\theta+j\pi)t\sqrt{-1}} \pmatrix
\sqrt{-1}a_j/2\\a_j/2\endpmatrix. \lb {3.12}\ee For
$\omega=e^{\sqrt{-1}\theta}$, $\theta\in [0,\pi)$, we define two
self-adjoint operators $A^{\omega}, B^{\omega}\in \mathcal
{L}(E^{\omega})$ by \bea  (A^{\omega}x,y)=\int^1_0\<-J\dot
x(t),y(t)\>dt,\;\;(B^{\omega}x,y)=\int^1_0\<B(t)x(t),y(t)\>dt
\nn\eea
 on $E^{\omega}$. Then $B^{\omega}$ is also compact.

\noindent{\bf Definition 3.2.} {\it We define the index function
 $$i_{\omega}^{L_0}(B)=I(A^{\omega}, \;\;A^{\omega}-B^{\omega}),\;\;\nu_{\omega}^{L_0}(B)=m^0(A^{\omega}-B^{\omega}),\;
  \forall\,\omega=e^{\sqrt{-1}\theta},\;\;\theta\in (0,\pi). $$
 }

 By the Floquet theory,
 we have $M^0(A^\om,B^\om)$ is isomorphic to the
 solution space  of the following linear Hamiltonian system
 $$\dot x(t)=JB(t)x(t) $$
 satisfying the following boundary condition
 $$x(0)\in L_0, \;\;x(1)\in e^{\theta J}L_0.$$
 If $m^0(A^\om,B^\om)>0$, there holds $$\gamma(1)L_0\cap e^{\theta J}L_0\neq
 \{0\}$$
 which is equivalent to  $$\omega^2=e^{2\theta\sqrt{-1}}\in
 \sigma\left([U(1)-\sqrt{-1}V(1)][U(1)+\sqrt{-1}V(1)]^{-1}\right).$$
 This claim follows from the fact that if $\gamma(1)L_0\cap e^{\theta J}L_0\neq
 \{0\}$, there exist $a, b\in \C^n\setminus \{0\}$ such that
 $$[U(1)+\sqrt{-1}V(1)]a=\omega^{-1}b,\;\;[U(1)-\sqrt{-1}V(1)]a=\omega b. $$
So we have  \be\nu_{\omega}^{L_0}(B)= \dim (\gamma(1)L_0\cap
e^{\theta J}L_0),\;\; \forall\,
\omega=e^{\sqrt{-1}\theta},\;\theta\in (0,\pi). \lb{ 3.14}\ee

\noindent{\bf Lemma 3.5.} {\it The index function
$i_{\omega}^{L_0}(B)$ is locally constant. For
$\omega_0=e^{\sqrt{-1}\theta_0},\;\theta_0\in (0,\pi)$ is a point of
discontinuity of $i_{\omega}^{L_0}(B)$, then
$\nu_{\omega_0}^{L_0}(B)>0$ and so $\dim (\gamma(1)L_0\cap
e^{\theta_0 J}L_0)>0$. Moreover there hold

\bea && |i_{\omega_0+}^{L_0}(B)-i_{\omega_0-}^{L_0}(B)|\le
\nu_{\omega_0}^{L_0}(B),\;\qquad
|i_{\omega_0+}^{L_0}(B)-i_{\omega_0}^{L_0}(B)|\le
\nu_{\omega_0}^{L_0}(B),\nn\\&&
|i_{\omega_0-}^{L_0}(B)-i_{\omega_0}^{L_0}(B)|\le
\nu_{\omega_0}^{L_0}(B),\;\qquad\;\;
|i_{L_0}(B)+n-i_{1+}^{L_0}(B)|\le \nu_{L_0}(B), \lb {3.15}\eea
 where $i_{\omega_0+}^{L_0}(B)$, $i_{\omega_0-}^{L_0}(B)$ are the
 limits on the right and left respectively of the index function $i_{\omega}^{L_0}(B)$
 at $\omega_0=e^{\sqrt{-1}\theta_0}$ as a function of $\theta$.}

\noindent{\bf Proof.} For $x(t)=e^{\theta tJ}u(t),
u(t)=\displaystyle\sum_{j\in\Z}e^{j\pi tJ}\pmatrix
0\\a_j\endpmatrix$, we have
$$((A^{\omega}-B^{\omega})x,x)=\int^1_{0}\<-J\dot u(t),u(t)\>dt+
\int^1_{0}\<(\theta-e^{-\theta tJ}B(t)e^{\theta tJ})u(t),u(t)\>dt.$$
 So we have
 $$((A^{\omega}-B^{\omega})x,x)=(q_{\omega}u,u) $$
 with
 $$(q_\om u,u)=
\int^1_{0}\<-J\dot u(t),u(t)\>dt+ \int^1_{0}\<(\theta-e^{-\theta
tJ}B(t)e^{\theta tJ})u(t),u(t)\>dt.
$$
 Since $\dim (\gamma(1)L_0\cap
e^{\theta J}L_0)>0$ at only finite (up to $n$) points $\theta\in
(0,\pi)$, for the point $\theta_0\in (0,\pi)$ such that
$\nu_{\omega_0}^{L_0}(B)=0$, then $\nu_{\omega}^{L_0}(B)=0$ for
$\omega=e^{\sqrt{-1}\theta}$, $\theta\in
(\theta_0-\delta,\theta_0+\delta)$, $\delta>0$ small enough. By
using the notations as in Lemma 3.3, we have
$$ (P_m^{\omega}(A^{\omega}-B^{\omega})P_m^{\omega}x,x)=(P_mq_{\omega}P_mu,u).$$
By Lemma 3.3, we have
$$m^0_d(P_m^{\omega}(A^{\omega}-B^{\omega})P_m^{\omega})=m^0(A^{\omega}-B^{\omega})=\nu_{\omega}^{L_0}(B)=0. $$
So by the continuity of the eigenvalue of a continuous family of
operators we have that
$$m^-_d(P_m^{\omega}(A^{\omega}-B^{\omega})P_m^{\omega}) $$
must be constant for $\omega=e^{\sqrt{-1}\theta}$, $\theta\in
(\theta_0-\delta,\theta_0+\delta)$. Since
$m^-_d(P_m^{\omega}A^{\omega}P_m^{\omega})$ is constant for
$\omega=e^{\sqrt{-1}\theta}$, $\theta\in
(\theta_0-\delta,\theta_0+\delta)$, we have $i_{\omega}^{L_0}(B)$ is
constant for $\omega=e^{\sqrt{-1}\theta}$, $\theta\in
(\theta_0-\delta,\theta_0+\delta)$.

The results in (\ref{3.15}) now follow from some standard arguments.
\hfill\hb

 By (\ref{c1}), Definition 3.2 and Lemma 3.5, we see that for
any $\omega_0=e^{\sqrt{-1}\theta_0},\;\theta_0\in (0,\pi)$, there
holds \be i^{L_0}_{\omega_0}(B)\ge
i_{L_0}(B)+n-\sum_{\omega=e^{\sqrt{-1}\theta},\; 0\le \theta\le
\theta_0}\nu^{L_0}_{\omega}(B). \lb {3.16}\ee
 We note that
 \be \sum_{\omega=e^{\sqrt{-1}\theta},\;
0\le \theta\le \theta_0}\nu^{L_0}_{\omega}(B)\le n. \lb{3.17}\ee
 So we have
 \be i_{L_0}(B)\le i^{L_0}_{\omega_0}(B)\le
i_{L_0}(B)+n. \lb {3.18}\ee

\setcounter{equation}{0}
\section{ Bott-type index formula for $L$-index} 
In this section, we establish the Bott-type  iteration formula for
the $L_j$-index theory with $j=0,1$. Without loss of generality, we
assume $\tau=1$. Suppose the continuous symplectic path $\gamma:
[0,1]\to \Sp(2n)$
 is the fundamental solution of the following linear Hamiltonian
 system
 \be \dot z(t)=J B(t)z(t),\quad t\in \R \lb{4.1}\ee
 with $B(t)$ satisfying $B(t+2)=B(t)$ and $B(1+t)N=NB(1-t))$ for $t\in \R$.  This implies
 $B(t)N=NB(-t)$ for $t\in \R$.  By the unique existence theorem of
 the linear differential equations, we get
 \be\gamma(1+t)=N\gamma(1-t)\gamma(1)^{-1}N\gamma(1), \gamma(2+t)=\gamma(t)\gamma(2). \lb{4.2}\ee
 For $j\in \N$, we define the $j$-times iteration path $\ga^j:[0,j]\to \Sp(2n)$ of $\gamma$ by
 $$\gamma^1(t)=\gamma(t), \;t\in [0,1], $$
 $$\gamma^2(t)=\left\{\begin{array}{l}  \gamma(t), \;t\in [0,1],\\
  N\gamma(2-t)\gamma(1)^{-1}N\gamma(1), \;t\in [1,2], \end{array}\right.$$
and in general, for $k\in\N$, we define
  \bea\gamma^{2k-1}(t)=\left\{\begin{array}{l} \gamma(t), \;t\in [0,1],\\
  N\gamma(2-t)\gamma(1)^{-1}N\gamma(1), \;t\in [1,2],\\\cdots\cdots\\
   N\gamma(2k-2-t)\gamma(1)^{-1}N\gamma(1)\gamma(2)^{2k-5}, \;t\in [2k-3,2k-2],\\
  \gamma(t-2k+2)\gamma(2)^{2k-4}, \;t\in [2k-2,2k-1],\end{array}\right.
  \lb{4.3}\eea
\bea\gamma^{2k}(t)=\left\{\begin{array}{l}\gamma(t), \;t\in [0,1],\\
 N\gamma(2-t)\gamma(1)^{-1}N\gamma(1), \;t\in [1,2],\\\cdots\cdots\\
   \gamma(t-2k+2)\gamma(2)^{2k-4}, \;t\in [2k-2,2k-1], \\
  N\gamma(2k-t)\gamma(1)^{-1}N\gamma(1)\gamma(2)^{2k-3}, \;t\in [2k-1,2k].\end{array}\right.\lb{4.4}\eea
For $\gamma\in \mathcal {P}_{\tau}(2n)$, we define \be\gamma^k(\tau
t)=\td{\gamma}^k(t)\;{\rm with}\; \td{\gamma}(t)=\gamma(\tau
t)\lb{uvw}.\ee
 For the $L_0$-index of the iteration path $\gamma^k$, we have
the following Bott-type formulas.

 \noindent{\bf Theorem 4.1.} {\it Suppose $\omega_k=e^{\pi
\sqrt{-1}/k}$. For odd $k$ we have
\bea i_{L_0}(\gamma^k)=i_{L_0}(\gamma^1)+\sum_{i=1}^{(k-1)/2}i_{\omega_k^{2i}}(\gamma^2),\nn\\
\nu_{L_0}(\gamma^k)=\nu_{L_0}(\gamma^1)+\sum_{i=1}^{(k-1)/2}\nu_{\omega_k^{2i}}(\gamma^2),\nn\eea
and for even $k$, we have \bea &&
i_{L_0}(\gamma^k)=i_{L_0}(\gamma^1)+i^{L_0}_{\omega_k^{k/2}}(\gamma^1)+
\sum_{i=1}^{k/2-1}i_{\omega_k^{2i}}(\gamma^2),\;
\nn\\
&&
\nu_{L_0}(\gamma^k)=\nu_{L_0}(\gamma^1)+\nu^{L_0}_{\omega_k^{k/2}}(\gamma^1)+
\sum_{i=1}^{k/2-1}\nu_{\omega_k^{2i}}(\gamma^2). \nn\eea}
 We note
that $\omega_k^{k/2}=\sqrt{-1}$.

Before proving Theorem 4.1, we give some notations and definitions.

We define the Hilbert space \bea E^k_{L_0}=\left\{x\in L^{2}([0,k],
\C^{2n})\,|\, x(t)=\sum_{j\in\Z}e^{jt\pi/kJ}\pmatrix
0\\a_j\endpmatrix,
\;a_j\in\C^{n},\;\|x\|^2:=\sum_{j\in\Z}(1+|j|)|a_j|^2<\infty\right\},
\nn\eea
 where we still denote  $L_0=\{0\}\times\C^n\subset \C^{2n}$ which is the Lagrangian
 subspace of the linear complex symplectic space
 $(\C^{2n},\omega_0)$.
 For $x\in E^k_{L_0}$, we can write
\bea x(t)&=&\sum_{j\in\Z}e^{jt\pi/kJ}\pmatrix 0\\a_j\endpmatrix
=\sum_{j\in \Z}\pmatrix
-\sin(jt\pi/k)a_j\\\cos(jt\pi/k)a_j\endpmatrix\nn\\&=&\sum_{j\in\Z}\left\{
e^{j\pi t\sqrt{-1}/k}\pmatrix
\sqrt{-1}a_j/2\\a_j/2\endpmatrix+e^{-j\pi t\sqrt{-1}/k}\pmatrix
-\sqrt{-1}a_j/2\\a_j/2\endpmatrix\right\}.\lb{4.5}\eea
 On $E^k_{L_0}$ we define two self-adjoint operators and a quadratical form by
 \be (A_kx,\,y)=\int_0^k\<-J\dot{x}(t),\,y(t)\>dt,\quad
 (B_kx,\,y)=\int_0^k\<B(t)x(t),y(t)\>dt, \lb{4.5'}\ee
 \be Q^k_{L_0}(x,y)=((A_k-B_k) x,y),\lb{4.6}\ee
 where in this section $\<\cdot,\cdot\>$ is the standard Hermitian inner product
 in $\C^{2n}$.

\noindent{\bf Lemma 4.1.} {\it $E_{L_0}^k$ has the following natural
decomposition

\be E^k_{L_0}=\bigoplus_{l=0}^{k-1}E_{L_0}^{\omega_k^l}, \lb{4.8}\ee
here we have extended the domain of functions in
 $E_{L_0}^{\omega_k^l}$ from $[0,1]$ to $[0,k]$ in the obvious
 way, i.e.,
$$E_{L_0}^{\omega_k^l}=\left\{x\in E^k_{L_0}\,|\,x(t)=e^{l\pi
tJ/k}\sum_{j\in\Z}e^{j\pi tJ}\pmatrix 0\\a_j\endpmatrix\right\}.$$
 }
\noindent{\bf Proof.}  Any element $x\in E^k_{L_0}$ can be written
as
 \bea  x(t)&=&\sum_{j\in\Z}\left\{
e^{j\pi t\sqrt{-1}/k}\pmatrix
\sqrt{-1}a_j/2\\a_j/2\endpmatrix+e^{-j\pi t\sqrt{-1}/k}\pmatrix
-\sqrt{-1}a_j/2\\a_j/2\endpmatrix\right\}\nn\\
&=&\sum_{l=0}^{k-1}\sum_{j\equiv l \,(mod k)}\left\{ e^{j\pi
t\sqrt{-1}/k}\pmatrix \sqrt{-1}a_j/2\\a_j/2\endpmatrix+e^{-j\pi
t\sqrt{-1}/k}\pmatrix
-\sqrt{-1}a_j/2\\a_j/2\endpmatrix\right\}\nn\\
&=&\sum_{l=0}^{k-1}\sum_{j\in\Z}\left\{e^{l\pi t\sqrt{-1}/k} e^{j\pi
t\sqrt{-1}}\pmatrix \sqrt{-1}b_j/2\\b_j/2\endpmatrix+e^{-l\pi
t\sqrt{-1}/k} e^{-j\pi t\sqrt{-1}}\pmatrix
-\sqrt{-1}b_j/2\\b_j/2\endpmatrix\right\}\nn\\&:=&
\xi_x(t)+N\xi_x(-t),\;\xi_x(t)=\sum_{l=0}^{k-1}\sum_{j\in\Z}e^{l\pi
t\sqrt{-1}/k} e^{j\pi t\sqrt{-1}}\pmatrix
\sqrt{-1}b_j/2\\b_j/2\endpmatrix,\lb{4.9}\eea
 where $b_j=a_{jk+l}$. By setting $\omega_k=e^{\pi \sqrt{-1}/k}$, and comparing (\ref{4.7}) and (\ref{4.9}), we
 obtain (\ref{4.8}).
\hfill\hb

Note that the natural decomposition (\ref{4.8}) is not orthogonal
under the quadratical form $Q_{L_0}^k$ defined in (\ref{4.6}). So
the type of the iteration formulas in Theorem 4.1 is somewhat
different from the original Bott formulas in \cite{Bott} of the
Morse index theory for closed geodesics and (\ref{1.21}) of
Maslov-type index theory for  periodic solutions of  Hamiltonian
systems and the Bott-type formulas in \cite{Ek}. This is also our
main difficulty in the proof of Theorem 4.1. However, after
recombining the terms in the decomposition in Lemma 4.1, we can
obtain an orthogonal decomposition under the quadratical form
$Q_{L_0}^k$.

For $1\le l<\frac {k}{2}$ and $l\in \N$, we set
$$E_{L_0}^{\omega_k,l}=E_{L_0}^{\omega_k^l}\oplus E_{L_0}^{\omega_k^{k-l}}.$$
So for odd $k$, we decompose $E_{L_0}^k$ as
$$E_{L_0}^k=E_{L_0}^1\oplus\bigoplus_{l=1}^{(k-1)/2}E_{L_0}^{\omega_k,l}, \eqno(C_{odd})$$
for even $k$, we decompose $E_{L_0}^k$ as
$$E_{L_0}^k=E_{L_0}^1\oplus E_{L_0}^{\omega_k^{k/2}}
\oplus\bigoplus_{l=1}^{\frac k2-1}E_{L_0}^{\omega_k,l}. \eqno(C_{even})$$

\noindent{\bf Lemma 4.2.} {\it The above two decompositions
($C_{odd}$) and ($C_{even}$) are orthogonal under the quadratical
form $Q_{L_0}^k$ for $k$ is odd and even respectively. Moreover, for
$x\in E_{L_0}^{\omega_k^i}$ and $y\in
 E_{L_0}^{\omega_k^j}$, $i,j\in\Z\cap[0,k-1]$, we have
  \bea &&(B_k
x,y)=\int^k_0 \<B(t)
x(t),y(t)\>\,dt=0, \;\;{ if}\; i\neq j,\;i+j\ne k,\lb{4.10}\\
&& (B_k x,y)=\int^k_0 \<B(t) x(t),y(t)\>\,dt\nn\\
&&\quad\quad\quad\;\;\,=k\int^1_0 \<B(t)
x(t),y(t)\>\,dt=k(B^{\om_k^i}x,y), \;\;{ if}\; i=j=
0,\frac {k}{2}, \lb{4.11}\\
 && (B_k x,y)=\int^k_0 \<B(t)
x(t),y(t)\>\,dt\nn\\
&&=k\left(\int^1_0 \<B(t) \xi_x(t),\xi_y(t)\>\,dt+\int^1_0 \<B(t)
N\xi_x(-t),N\xi_y(-t)\>\,dt\right), \;{ if}\; i=j\neq 0,\frac
{k}{2}, \lb{4.12}\eea \bea && (B_kx,y)=k\left(\int^1_0 \<B(t)
N\xi_x(-t),\xi_y(t)\>\,dt\right.\nn\\
&&\qquad\qquad\;\;\quad\quad\left.+\int^1_0 \<B(t)
\xi_x(t),N\xi_y(-t)\>\,dt\right),\;\;{ if}\; i\ne j,\; i+j=k,\lb{4.13}\\
&& (A_k x,y)=\int^k_0 \<-J\dot x(t),y(t)\>\,dt=0, \;\;{ if}\; i\neq j,\lb{4.13'}\\
&&  (A_k x,y)=\int^k_0 \<-J\dot x(t),y(t)\>\,dt=k\int^1_0 \<-J\dot
x(t),y(t)\>\,dt=k(A^{\om_k^i}x,y), \;\;{ if}\; i=j, \lb{4.14}\eea
 where the operators $A^{\omega}$, $B^{\omega}$ are defined
 in Section 3.}

\noindent{\bf Proof.} We first prove the formulas
(\ref{4.10})-(\ref{4.14}). It is easy to see that, we only need to
prove them in the case
 \bea &&
x(t)=e^{it\pi\sqrt{-1}/k}e^{pt\pi
\sqrt{-1}}\alpha_p+e^{-it\pi\sqrt{-1}/k}e^{-pt\pi
\sqrt{-1}}N\alpha_p,\nn\\
&& y(t)=e^{jt\pi\sqrt{-1}/k}e^{mt\pi
\sqrt{-1}}\alpha_m+e^{-jt\pi\sqrt{-1}/k}e^{-mt\pi
\sqrt{-1}}N\alpha_m,\nn\\&& \alpha_s=\pmatrix \sqrt{-1}
a_s\\a_s\endpmatrix,\nn\eea for any integers $p$ and $m$.

In this case,
\bea (B_kx,y)&=&\int^k_0\<B(t)\alpha_p,\;e^{(j-i)t\pi\sqrt{-1}/k}e^{(m-p)t\pi\sqrt{-1}}\alpha_m\>\,dt\nn\\
&\,&\;+\int^k_0\<B(t)\alpha_p,\;e^{-(j+i)t\pi\sqrt{-1}/k}e^{-(m+p)t\pi\sqrt{-1}}N\alpha_m\>\,dt\nn\\
&\,&\;+\int^k_0\<B(t)N\alpha_p,\;e^{(j+i)t\pi\sqrt{-1}/k}e^{(m+p)t\pi\sqrt{-1}}\alpha_m\>\,dt\nn\\
&\,&\;+\int^k_0\<B(t)N\alpha_p,\;e^{(i-j)t\pi\sqrt{-1}/k}e^{(p-m)t\pi\sqrt{-1}}N\alpha_m\>\,dt\nn\\
&=&\sum_{s=1}^k\int^s_{s-1}\<B(t)\alpha_p,\;e^{(j-i)t\pi\sqrt{-1}/k}e^{(m-p)t\pi\sqrt{-1}}\alpha_m\>\,dt\nn\\
&\,&\;+\sum_{s=1}^k\int^s_{s-1}\<B(t)\alpha_p,\;e^{-(j+i)t\pi\sqrt{-1}/k}e^{-(m+p)t\pi\sqrt{-1}}N\alpha_m
\>\,dt\nn\\
&\,&\;+\sum_{s=1}^k\int^s_{s-1}\<B(t)N\alpha_p,\;e^{(j+i)t\pi\sqrt{-1}/k}e^{(m+p)t\pi\sqrt{-1}}\alpha_m)\>\,dt\nn\\
&\,&\;+\sum_{s=1}^k\int^s_{s-1}\<B(t)N\alpha_p,\;e^{(i-j)t\pi\sqrt{-1}/k}e^{(p-m)t\pi\sqrt{-1}}N\alpha_m\>\,dt\nn\\
&:=&I_1+I_2+I_3+I_4.\nn\eea
 By using the relations $B(1+t)N=NB(1-t)$
and $B(t)N=NB(-t)$, we have
$$$$
\bea && \int^{s+1}_{s}\<B(t)\alpha_p,\;e^{(j-i)t\pi\sqrt{-1}/k}e^{(m-p)t\pi\sqrt{-1}}\alpha_m\>\,dt\nn\\
&=&\int^{s}_{s-1}\<B(1+t)\alpha_p,\;e^{(j-i)(1+t)\pi\sqrt{-1}/k}e^{(m-p)(1+t)\pi\sqrt{-1}}\alpha_m\>\,dt\nn\\
&=&\int^{s}_{s-1}\<NB(1-t)N\alpha_p,\;e^{(j-i)(1+t)\pi\sqrt{-1}/k}e^{(m-p)(1+t)\pi\sqrt{-1}}\alpha_m\>\,dt\nn\\
&=&\int^{s}_{s-1}\<B(t-1)\alpha_p,\;e^{(j-i)(1+t)\pi\sqrt{-1}/k}e^{(m-p)(1+t)\pi\sqrt{-1}}\alpha_m\>\,dt\nn\\
&=&\int^{s-1}_{s-2}\<B(t)\alpha_p,\;e^{(j-i)(2+t)\pi\sqrt{-1}/k}e^{(m-p)(2+t)\pi\sqrt{-1}}\alpha_m\>\,dt\nn\\
&=&e^{2(i-j)\pi\sqrt{-1}/k}\int^{s-1}_{s-2}\<B(t)\alpha_p,\;
e^{(j-i)t\pi\sqrt{-1}/k}e^{(m-p)t\pi\sqrt{-1}}\alpha_m\>\,dt.\nn\eea
Similarly, we have
\bea &&\int^{s+1}_{s}\<B(t)\alpha_p,\;e^{-(j+i)t\pi\sqrt{-1}/k}e^{-(m+p)t\pi\sqrt{-1}}N\alpha_m\>\,dt\nn\\
&=&e^{2(j+i)\pi\sqrt{-1}/k}\int^{s-1}_{s-2}
\<B(t)\alpha_p,\;e^{-(j+i)t\pi\sqrt{-1}/k}e^{-(m+p)t\pi\sqrt{-1}}N\alpha_m\>\,dt.\nn\\
 &&\int^{s+1}_{s}\<B(t)N\alpha_p,\;e^{(j+i)t\pi\sqrt{-1}/k}e^{(m+p)t\pi\sqrt{-1}}\alpha_m\>\,dt\nn\\
&=&e^{-2(j+i)\pi\sqrt{-1}/k}\int^{s-1}_{s-2}
\<B(t)N\alpha_p,\;e^{-(j+i)t\pi\sqrt{-1}/k}e^{-(m+p)t\pi\sqrt{-1}}\alpha_m\>\,dt.\nn\\
 &&\int^{s+1}_{s}\<B(t)N\alpha_p,\;e^{(i-j)t\pi\sqrt{-1}/k}e^{(p-m)t\pi\sqrt{-1}}N\alpha_m\>\,dt\nn\\
&=&e^{2(j-i)\pi\sqrt{-1}/k}\int^{s-1}_{s-2}
\<B(t)N\alpha_p,\;e^{(i-j)t\pi\sqrt{-1}/k}e^{(p-m)t\pi\sqrt{-1}}N\alpha_m\>\,dt.\nn\\
 &&\int^{2}_{1}\<B(t)\alpha_p,\;e^{(j-i)t\pi\sqrt{-1}/k}e^{(m-p)t\pi\sqrt{-1}}\alpha_m\>\,dt\nn\\
&=&e^{2(i-j)\pi\sqrt{-1}/k}\int^{1}_{0}
\<B(t)N\alpha_p,\;e^{(i-j)t\pi\sqrt{-1}/k}e^{(p-m)t\pi\sqrt{-1}}N\alpha_m\>\,dt.\nn\\
 &&\int^{2}_{1}\<B(t)\alpha_p,\;e^{-(j+i)t\pi\sqrt{-1}/k}e^{-(m+p)t\pi\sqrt{-1}}N\alpha_m\>\,dt\nn\\
&=&e^{2(j+i)\pi\sqrt{-1}/k}\int^{1}_{0}
\<B(t)N\alpha_p,\;e^{(j+i)t\pi\sqrt{-1}/k}e^{(m+p)t\pi\sqrt{-1}}\alpha_m\>\,dt.\nn\\
&&\int^{2}_{1}\<B(t)N\alpha_p,\;e^{(j+i)t\pi\sqrt{-1}/k}e^{(m+p)t\pi\sqrt{-1}}\alpha_m\>\,dt\nn\\
&=&e^{-2(j+i)\pi\sqrt{-1}/k}\int^{1}_{0}
\<B(t)\alpha_p,\;e^{-(j+i)t\pi\sqrt{-1}/k}e^{-(m+p)t\pi\sqrt{-1}}N\alpha_m\>\,dt.\nn\\
&&\int^{2}_{1}\<B(t)N\alpha_p,\;e^{(i-j)t\pi\sqrt{-1}/k}e^{(p-m)t\pi\sqrt{-1}}N\alpha_m\>\,dt\nn\\
&=&e^{2(j-i)\pi\sqrt{-1}/k}\int^{1}_{0}
\<B(t)\alpha_p,\;e^{(j-i)t\pi\sqrt{-1}/k}e^{(m-p)t\pi\sqrt{-1}}\alpha_m\>\,dt.\nn\eea
From these observations, we find that
$$\;I_2+I_3=0, \;{\rm if }\;i+j\neq 0,k$$
and
$$I_1+I_4=0, \;{\rm if }\;i\neq j$$
which yield (\ref{4.10}).
 In fact, by setting $\mu=e^{2(i-j)\pi \sqrt{-1}/k}$, then
$\mu^k=1$, for $k=2q$ with $q\in \N$, we have
\bea
I_1&=&(1+\mu+\cdots+\mu^{q-1})\int^1_0\<B(t)\alpha_p,\;e^{(j-i)t\pi\sqrt{-1}/k}
e^{(m-p)t\pi\sqrt{-1}}\alpha_m\>\,dt
\nn\\&&+(\mu+\cdots+\mu^q)\int^1_0
\<B(t)N\alpha_p,\;e^{(i-j)t\pi\sqrt{-1}/k}e^{(p-m)t\pi\sqrt{-1}}N\alpha_m\>\,dt.\nn\eea
\bea
I_4&=&(\mu^{-1}+\cdots+\mu^{-q})\int^1_0\<B(t)\alpha_p,\;e^{(j-i)t\pi\sqrt{-1}/k}e^{(m-p)t\pi\sqrt{-1}}\alpha_m\>\,dt
\nn\\&&+(1+\mu^{-1}+\cdots+\mu^{-q+1})\int^1_0
\<B(t)N\alpha_p,\;e^{(i-j)t\pi\sqrt{-1}/k}e^{(p-m)t\pi\sqrt{-1}}N\alpha_m\>\,dt.\nn\eea
Noting
$$\mu^{-1}+\cdots+\mu^{-q}+1+\mu+\cdots+\mu^{q-1}=\frac{\mu^{-q}(1-\mu^{2q})}{1-\mu}=0$$
and
$$\mu+\cdots+\mu^q+1+\mu^{-1}+\cdots+\mu^{-q+1}=\frac{\mu^{-q+1}(1-\mu^{2q})}{1-\mu}=0,$$
we have $I_1+I_4=0$ provided $i-j\neq 0$. For $k=2q-1$ with $q\in
\N$, in the similar way we also have $I_1+I_4=0$ provided $i-j\neq
0$. That $I_2+I_3=0$ provided $i+j\neq 0,k$ is proved in the same
way.

For the case $i=j=0$ and the case $i=j=\frac{k}{2}$ if $k$ is even,
from the above observation we have
$$\int^k_0\<B(t)x(t),\;y(t)\>dt=k\int^1_0 \<B(t)x(t),\;y(t)\>dt$$
which yields (\ref{4.11}).

For the cases $i=j\neq 0,\frac{k}{2}$, we have $I_2+I_3=0$ and
\bea (B_kx,y)&=&I_1+I_4\nn\\
&=&k\left(\int^1_0\<B(t)\alpha_p,\;e^{(j-i)t\pi\sqrt{-1}/k}e^{(m-l)t\pi\sqrt{-1}}\alpha_m\>\,dt\right.\nn\\
&&\;\;\;\;\;+\left.\int^1_0\<B(t)N\alpha_p,\;e^{(i-j)t\pi\sqrt{-1}/k}e^{(l-m)t\pi\sqrt{-1}}N\alpha_m\>\,dt\right)\nn\\
&=&k\left(\int^1_0\<B(t)\xi_x(t),\;\xi_y(t)\>\,dt+\int^1_0\<B(t)N\xi(-t),\;N\eta(-t)\>\,dt\right),\lb{4.15}\eea
where for $x,y\in E^{\omega_k^i}_{L_0}$, $\xi_x$ and $\xi_y$ are
defined in as in (\ref{4.9}). So (\ref{4.12}) holds from
(\ref{4.15}). The claim (\ref{4.13}) is proved by the same way. By
direct computation we have (\ref{4.13'}) and (\ref{4.14}), moreover
$$(A_kx,y)=k\left(\int^1_0\<-J\frac {d}{dt}\xi_x(t),\;\xi_y(t)\>\,dt+\int^1_0\<-J\frac{d}{dt}N\xi_x(-t),
\;N\xi_y(-t)\>\,dt\right),\; {\rm if }\; i=j.$$ The orthogonality
statement in Lemma 4.2 follows from (\ref{4.10}) and (\ref{4.13'}).

\hfill\hb

\noindent{\bf Proof of Theorem 4.1.} Let $1\le l<\frac{k}{2}$, $l\in
\N$.
 For $x\in
E_{L_0}^{\omega_k^{l}}$,
$$ x(t)=\sum_{j\in\Z}e^{l\pi\sqrt{-1}t/k}e^{j\pi\sqrt{-1}t}\pmatrix \sqrt{-1}\alpha_j\\\alpha_j\endpmatrix+
e^{-l\pi\sqrt{-1}t/k}e^{-j\pi\sqrt{-1}t}\pmatrix
-\sqrt{-1}\alpha_j\\\alpha_j\endpmatrix.$$
 For $y\in E_{L_0}^{\omega_k^{k-l}}$,
 $$ y(t)=\sum_{j\in\Z}e^{-l\pi\sqrt{-1}t/k}e^{-j\pi\sqrt{-1}t}\pmatrix \sqrt{-1}\beta_j\\\beta_j\endpmatrix+
e^{l\pi\sqrt{-1}t/k}e^{j\pi\sqrt{-1}t}\pmatrix
-\sqrt{-1}\beta_j\\\beta_j\endpmatrix.$$
 Thus for $z=x+y\in E_{L_0}^{\omega_k,l}$ with $x\in
 E_{L_0}^{\omega_k^{l}}$ and $y\in E_{L_0}^{\omega_k^{k-l}}$,

\bea  z(t)&=&
\sum_{j\in\Z}e^{l\pi\sqrt{-1}t/k}e^{j\pi\sqrt{-1}t}\pmatrix
\sqrt{-1}\alpha_j\\\alpha_j\endpmatrix+
e^{-l\pi\sqrt{-1}t/k}e^{-j\pi\sqrt{-1}t}\pmatrix
-\sqrt{-1}\alpha_j\\\alpha_j\endpmatrix\nn\\
&\;&+e^{-l\pi\sqrt{-1}t/k}e^{-j\pi\sqrt{-1}t}\pmatrix
\sqrt{-1}\beta_j\\\beta_j\endpmatrix+
e^{l\pi\sqrt{-1}t/k}e^{j\pi\sqrt{-1}t}\pmatrix
-\sqrt{-1}\beta_j\\\beta_j\endpmatrix\nn\\
& =&\xi_x(t)+N\xi_x(-t)+\xi_y(-t)+N\xi_y(t).\nn\eea
 So for $z=x+y\in E_{L_0}^{\omega_k,l}$ with $x\in
 E_{L_0}^{\omega_k^{l}}$ and $y\in E_{L_0}^{\omega_k^{k-l}}$, we have
 \bea
 (B_kz,z)&=&(B_kx,x)+(B_ky,y)+(B_kx,y)+(B_ky,x)\nn\\
 &=&k\left(\int_0^1\<B(t)\xi_x(t),\;\xi_x(t)\>dt+\int_0^1\<B(t)\xi_x(t),\;N\xi_y(t)\>dt+\right.\nn\\
 &&+\int_0^1\<B(t)N\xi_x(-t),\;N\xi_x(-t)\>dt+\int_0^1\<B(t)N\xi_x(-t),\;\xi_y(-t)\>dt+\nn\\
 &&+\int_0^1\<B(t)\xi_y(-t),\;\xi_y(-t)\>dt+\int_0^1\<B(t)\xi_y(-t),\;N\xi_x(-t)\>dt+\nn\\
 &&+\left.\int_0^1\<B(t)N\xi_y(t),\;N\xi_y(t)\>dt+\int_0^1\<B(t)N\xi_y(t),\;\xi_x(t)\>dt\right)\nn\\
 &=& k\int^1_{-1}\<B(t)(\xi_x(t)+N\xi_y(t)),\; \xi_x(t)+N\xi_y(t)\>dt\nn\\
 &=& k\int^2_{0}\<B(t)(\xi_x(t)+N\xi_y(t)),\; \xi_x(t)+N\xi_y(t)\>dt,\nn\eea
where in the second equality we have used (\ref{4.12}) and
(\ref{4.13}).

 We note that
 \bea u(t)&=&\xi_x(t)+N\xi_y(t)=\displaystyle\sum_{j\in\Z}e^{l\pi\sqrt{-1}t/k}e^{j\pi\sqrt{-1}t}
 \pmatrix
 \sqrt{-1}(\alpha_j-\beta_j)\\(\alpha_j+\beta_j)\endpmatrix\nn\\&=&\sum_{j\in\Z}e^{l\pi\sqrt{-1}t/k}e^{j\pi\sqrt{-1}t}u_j,\;\;
 u_j\in \C^{2n}.\nn\eea We set
 $$E_{\omega_k^{2l} }=\left\{u\in L^{2}([0,2],\C^{2n})\,|\,u(t)=e^{l\pi\sqrt{-1}t/k}
 \sum_{j\in\Z}e^{j\pi\sqrt{-1}t}u_j, \;\|u\|^2:=\sum_{j\in\Z}(1+|j|)|u_j|^2<+\infty\right\}. $$
 We define self-adjoint operators  on $E_{\omega_k^{2l} }$ by
 $$(A_{\omega_k^{2l} }u,v)=\int^2_0\<-J\dot u(t),\;v(t)\>dt,\;(B_{\omega_k^{2l} }u,v)=\int^2_0\<B(t) u(t),
 \;v(t)\>dt$$
 and a quadratic form
 $$Q_{\omega_k^{2l} }(u)=((A_{\omega_k^{2l} }-B_{\omega_k^{2l} })u,u), \;u\in E_{\omega_k^{2l} }. $$
 Here $Q_{\omega }$ is just the quadratic form $f_{\omega}$ defined
 on p$_{133}$ of
 \cite{Long1}. In order to complete the proof of Theorem 4.1, we
 need the following result.

\noindent{\bf Lemma 4.3.} {\it For a symmetric 2-periodic matrix
function $B$ and $\om\in \U\setminus\{1\}$, there hold
\bea
&&I(A_{{\omega }},A_{\omega }-B_{\omega })=i_{\omega
}(\ga^2),\lb{b1}\\
      &&m^0(A_{\omega}-B_{\omega })=\nu_{\omega }(\ga^2).\lb{b2}\eea}
\noindent{\bf Proof.}  In fact, (\ref{b1}) follows directly from
Definition 2.3 and Corollary 2.1 of \cite{LZ1} and Lemma 3.4,
(\ref{b2}) follows from the Floquet theory. We note also that
(\ref{b1}) is the eventual form of the Galerkin approximation
formula. We can also prove it step by step as the proof of Theorem
3.1 of \cite{Liu0} by using the saddle point reduction formula  in
Theorem 6.1.1 of \cite{Long1}. \hfill\hb

\noindent{\it Continue the proof of Theorem 4.1. } By Lemma 4.3, we
have
 \bea I(A_{{\omega_k^{2l} }},A_{\omega_k^{2l} }-B_{\omega_k^{2l} })=i_{\omega_k^{2l} }(\ga^2),\;\;
      m^0(A_{\omega_k^{2l} }-B_{\omega_k^{2l} })=\nu_{\omega_k^{2l} }(\ga^2),\;\,1\le l< \frac{k}{2},\;\,l\in \N.\lb{4.16}\eea
By Definition 3.2, we have
    \bea I(A^{\sqrt{-1}},A^{\sqrt{-1}}-B^{\sqrt{-1}})=i_{\sqrt{-1}}^{L_0}(\ga),\;\;\;\;m^0(A^{\sqrt{-1}}-B^{\sqrt{-1}})=
    \nu_{\sqrt{-1}}^{L_0}(\ga).\lb{4.17}\eea
 By (\ref{c1}) we have
     \be I(A^1,A^1-B^1)=i_{L_0}(\ga)+n,\;\;\;\;m^0(A^1-B^1)=\nu_{L_0}(\ga),\lb{4.18}\ee
and  \be
I(A_k,A_k-B_k)=i_{L_0}(\ga^k)+n,\;\;\;\;m^0(A_k-B_k)=\nu_{L_0}(\ga^k)
.\lb{4.19}\ee By (\ref{4.11}), (\ref{4.14}), Lemma 3.3, Definition
3.1 and Lemma 4.2, for  odd $k$, sum the first equality in
(\ref{4.16}) for $l=1,2,\cdots,\frac{k-1}{2}$ and the first equality
of (\ref{4.18}) correspondingly. By comparing with the first
equality of (\ref{4.19}) we have
 \be
 i_{L_0}(\ga^k)=i_{L_0}(\ga)+\sum_{l=1}^{\frac{k-1}{2}}i_{\om_k^{2l}}(\ga^2),\lb{4.20}\ee
and for  even $k$, sum the first equality in (\ref{4.16}) for
$l=1,2,\cdots,\frac{k}{2}-1$ and the first equalities of
(\ref{4.17})-(\ref{4.18}) correspondingly. By comparing with the
first equality of (\ref{4.19}) we have
 \be
 i_{L_0}(\ga^k)=i_{L_0}(\ga)+i_{\sqrt{-1}}^{L_0}(\ga)+\sum_{l=1}^{\frac{k}{2}-1}i_{\om_k^{2l}}(\ga^2).\lb{4.21}\ee
Similarly we have \bea
&&\nu_{L_0}(\ga^k)=\nu_{L_0}(\ga)+\sum_{l=1}^{\frac{k-1}{2}}\nu_{\om_k^{2l}}(\ga^2),\quad {\rm if \; k\; is\; odd},\lb{4.22}\\
&&\nu_{L_0}(\ga^k)=\nu_{L_0}(\ga)+\nu_{\sqrt{-1}}^{L_0}(\ga)+\sum_{l=1}^{\frac{k}{2}-1}\nu_{\om_k^{2l}}(\ga^2),\quad
{\rm if \; k\; is\; even}.\lb{4.23}\eea
Then Theorem 4.1 holds from
(\ref{4.20})-(\ref{4.23}) and the fact that $\om_k^{k/2}=\sqrt{-1}$.
\hfill\hb

 From the formulas in Theorem 4.1, we note that
$$i_{L_0}(\gamma^2)=i_{L_0}(\gamma^1)+i^{L_0}_{\sqrt{-1}}(\gamma^1),\;\;
\nu_{L_0}(\gamma^2)=\nu_{L_0}(\gamma^1)+\nu^{L_0}_{\sqrt{-1}}(\gamma^1).$$
It implies (\ref{1.20}).

\noindent{\bf Definition 4.1.} {\it The mean $L_0$-index of $\gamma$
is defined
 by
 $$\hat {i}_{L_0}(\gamma)=\lim_{k\to +\infty}\frac{i_{L_0}(\gamma^k)}{k}.
 $$}
By definitions of $\hat{i}_{L_0}(\ga)$ and $\hat{i}(\ga^2)$(cf.
\cite{Long1} for example), the following result is obvious.

 \noindent{\bf Proposition 4.1.} {\it The mean $L_0$-index of $\gamma$ is
well defined, and
\bea\hat {i}_{L_0}(\gamma)=\frac
{1}{2\pi}\int^{\pi}_0i_B(e^{\sqrt{-1}\theta})d\theta=\frac{\hat
{i}(\gamma^2)}{2}, \lb{4.24}\eea  }
 here we have written
$i_B(\omega)=i_{\omega}(B)=i_{\omega}(\gamma_B)$.

 For $L_1=\R^n\times\{0\}$, we have the $L_1$-index theory
established in \cite{Liu2}. Similarly as in Definition 3.2, for
$\omega=e^{\theta\sqrt{-1}},\;\theta\in(0,\pi)$, we define
 $$E^{\omega}_{L_1}=\left\{x\in
L^2([0,1],\C^{2n})\,|\,x(t)=e^{\theta t
J}\displaystyle\sum_{j\in\Z}e^{j\pi tJ}\pmatrix a_j\\0\endpmatrix,
\;a_j\in\C^n,\;
\|x\|:=\sum_{j\in\Z}(1+|j|)|a_j|^2<+\infty\right\}.$$
 In
$E^{\omega}_{L_1}$ we define two  operators $A^{\omega}_{L_1}$ and
 $B^{\omega}_{L_1}$  by the same way  as the definitions of operators $A^{\omega}$ and
 $B^{\omega}$  in the section 3, but the
 domain is $E^{\omega}_{L_1}$. We define
$$i^{L_1}_{\omega}(B)=I(A^{\omega}_{L_1},A^{\omega}_{L_1}-B^{\omega}_{L_1}),\;\nu^{L_1}_\omega(B)=m^0(A^{\omega}_{L_1}-B^{\omega}_{L_1})). $$
\noindent{\bf Theorem 4.2.} {\it Suppose $\omega_k=e^{\pi
\sqrt{-1}/k}$. For odd $k$ we have \bea &&
i_{L_1}(\gamma^k)=i_{L_1}(\gamma^1)+\sum_{i=1}^{\frac{k-1}{2}}i_{\omega_k^{2i}}(\gamma^2),\nn\\
&&\nu_{L_1}(\gamma^k)=\nu_{L_1}(\gamma^1)+\sum_{i=1}^{\frac{k-1}{2}}\nu_{\omega_k^{2i}}(\gamma^2).\eea
 For even $k$, we have
\bea  &&
i_{L_1}(\gamma^k)=i_{L_1}(\gamma^1)+i^{L_1}_{\omega_k^{k/2}}(\gamma^1)+
\sum_{i=1}^{k/2-1}i_{\omega_k^{2i}}(\gamma^2),\;
\nn\\
&&
\nu_{L_1}(\gamma^k)=\nu_{L_1}(\gamma^1)+\nu^{L_1}_{\omega_k^{k/2}}(\gamma^1)+
\sum_{i=1}^{k/2-1}\nu_{\omega_k^{2i}}(\gamma^2). \nn\eea }
\noindent{\bf Proof.} The proof is almost the same as  that of
Theorem 4.1. The only thing different from that is the matrix $N$
should be replaced by $N_1=-N$. \hfill\hb

It is easy to see that
$i(\gamma^2)=i_{L_0}(\gamma^1)+i_{L_1}(\gamma^1)+n$, see Proposition
C of \cite{LZZ} for a proof, we remind that
$\mu_1(\gamma)=i_{L_0}(\gamma)+n$ and
$\mu_2(\gamma)=i_{L_1}(\gamma)+n$ (see (\ref{6.9}) below). So by the
Bott-type formula (see \cite{Long0}) for the $\omega$-index of
$\gamma^2$ at $\omega=-1$, we have
$$i_{-1}(\gamma^2)=i^{L_0}_{\sqrt{-1}}(\gamma^1)+i^{L_1}_{\sqrt{-1}}(\gamma^1), $$
$$\nu_{-1}(\gamma^2)=\nu^{L_0}_{\sqrt{-1}}(\gamma^1)+\nu^{L_1}_{\sqrt{-1}}(\gamma^1). $$

We now give a direct  proof of this result.

\noindent{\bf Proposition 4.2.} {\it There hold \bea
&& i(\gamma^2)=i_{L_0}(\gamma^1)+i_{L_1}(\gamma^1)+n, \lb{4.25}\\
&& \nu_1(\gamma^2)=\nu_{L_0}(\gamma^1)+\nu_{L_1}(\gamma^1), \lb{4.26}\\
&& i_{-1}(\gamma^2)=i^{L_0}_{\sqrt{-1}}(\gamma^1)+i^{L_1}_{\sqrt{-1}}(\gamma^1),\lb{4.27}\\
&&
\nu_{-1}(\gamma^2)=\nu^{L_0}_{\sqrt{-1}}(\gamma^1)+\nu^{L_1}_{\sqrt{-1}}(\gamma^1).
\lb{4.28}\eea } \noindent{\bf Proof.} Set $E_1=W^{1/2,2}(S^1,
\C^{2n})$ with $S^1=\R/(2\Z)$. We note that $E_{\omega}=e^{J\theta
t}E_1$ for $\omega=e^{2\theta\sqrt{-1}}$. For any $z\in E_1$, we
have
$$z(t)=\sum_{j\in\Z}e^{jt\pi J}c_j=\sum_{j\in\Z}e^{jt\pi J}\pmatrix 0\\a_j\endpmatrix
+\sum_{j\in\Z}e^{jt\pi J}\pmatrix b_j\\0\endpmatrix, \;c_j\in
\C^{2n},\;a_j,\;b_j\in \C^{n}.$$
 So we have $E_{\omega}=E_{L_0}^{\omega}\oplus E_{L_1}^{\omega}$.
 For $x\in E_{L_0}^{\omega}$ and $y\in E_{L_1}^{\omega}$, we can
 write
 \bea  x(t)&=&e^{J\theta t} \sum_{j\in\Z}e^{jt\pi J}\pmatrix 0\\a_j\endpmatrix:=e^{J\theta t}x_0(t),\nn\\
 y(t)&=&e^{J\theta t} \sum_{j\in\Z}e^{jt\pi J}\pmatrix b_j\\0\endpmatrix:=e^{J\theta
 t}y_0(t).\nn\eea
 By setting $\tilde {B}(t)=e^{-J\theta t}B(t)e^{J\theta t}$, we get
 $$\int^2_0\<B(t)x(t),y(t)\>dt=\int^2_0\<\tilde {B}(t)x_0(t),y_0(t)\>dt. $$
 In the cases of $\theta=0,\frac{\pi}{2}$, we have $\tilde {B}(t+2)=\tilde
 {B}(t)$ and $\tilde B(1+t)=N\tilde B(1-t)N$. As in (\ref{3.12}), we write
 $x_0(t)=\xi(t)+N\xi(-t)$ and  $y_0(t)=\eta(t)-N\eta(-t)$ with
 $$\xi(t)=\sum_{j\in \Z} e^{j\pi t\sqrt{-1}}\pmatrix \sqrt{-1}a_j\\a_j\endpmatrix,\;
 \eta(t)=\sum_{j\in \Z} e^{j\pi t\sqrt{-1}}\pmatrix b_j\\-\sqrt{-1}b_j\endpmatrix.$$

 \bea &&\int^2_1\<\tilde B(t)x_0(t),\;y_0(t)\>dt=\int^2_1\<\tilde B(t)(\xi(t)+N\xi(-t)),\;\eta(t)-N\eta(-t)\>dt\nn\\
 &&=
  \sum_{j,l\in\Z}\int^1_0\left\<\tilde B(1+t)\left(e^{j\pi(t+1)\sqrt{-1}}\pmatrix
  \sqrt{-1}a_j\\a_j\endpmatrix+
  e^{-j\pi(t+1)\sqrt{-1}}\pmatrix -\sqrt{-1}a_j\\a_j\endpmatrix\right)\right.,\;\nn\\
  && \quad\quad\quad\quad\quad\quad
  \left. e^{l\pi(t+1)\sqrt{-1}}\pmatrix b_j\\-\sqrt{-1}b_j
  \endpmatrix+e^{-l\pi(t+1)\sqrt{-1}}\pmatrix b_j\\\sqrt{-1}b_j\endpmatrix\right\>dt\nn\\
  &&=
   \sum_{j,l\in\Z}(-1)^{j+l}\int^1_0\<N\tilde B(1-t)N(\xi(t)+N\xi(-t)),\;\eta(t)-N\eta(-t)\>dt\nn\\
   &&=
    \sum_{j,l\in\Z}(-1)^{j+l}\int^1_0\<N\tilde B(t)N(\xi(1-t)+N\xi(t-1)),\;\eta(1-t)-N\eta(t-1)\>dt\nn\\
    &&=
     \sum_{j,l\in\Z}(-1)^{2(j+l)}\int^1_0\<\tilde B(t)(N\xi(-t)+\xi(t)),\;-\eta(t)+N\eta(-t)\>dt\nn\\
     &&=
     -\int^1_0\<\tilde B(t)(\xi(t)+N\xi(-t)),\;\eta(t)-N\eta(-t)\>dt=-\int^1_0\<\tilde
     B(t)x_0(t),\;y_0(t)\>dt.\nn\eea
     It implies that
    \be\int^2_0\<\tilde B(t)x_0(t),\;y_0(t)\>dt=0. \lb{4a}\ee
It is easy to see that
 \be\int^2_0\<-J\dot x(t),\;y(t)\>dt=0.\lb{4b} \ee
  By
defining
$$Q_{\omega}(x,y)=\int^2_0\<-J\dot x(t),\;y(t)\>dt-\int^2_0\<B(t)x(t),\;y(t)\>dt, \;x,\;y\in E_{\omega}, $$
(\ref{4a}) and (\ref{4b}) imply that the decomposition
$E_{\omega}=E_{L_0}^{\omega}\oplus E_{L_1}^{\omega}$ is
$Q_{\omega}$-orthogonal in the cases $\theta=0, \frac{\pi}{2}$. So
we get the formulas  (\ref{4.25})-(\ref{4.28}) by the similar
argument in the proof of Theorem 4.1. \hfill\hb

\setcounter{equation}{0}
\section {Proof of Theorems 1.4 and 1.5} 

\noindent{\bf Proof of Theorem 1.4.} By the definition of the
splitting number, we have
$$i_{\omega_0}(\gamma^2)=i(\gamma^2)+\sum_{0\le \theta<\theta_0}S^+_M(e^{\sqrt{-1}\theta})-
\sum_{0< \theta\le \theta_0}S^-_M(e^{\sqrt{-1}\theta}), $$ where
$\omega_0=e^{\sqrt{-1}\theta_0}$. So for $k\in 2\N-1$, let
$m=\frac{k-1}{2}$, we have
 \bea && \sum_{i=1}^m i_{\omega_k^{2i}}(\gamma^2)=m
 i(\gamma^2)+\sum_{i=1}^m\left(\sum_{0\le \theta<\frac{2i\pi}{k}}S^+_M(e^{\sqrt{-1}\theta})-
 \sum_{0< \theta\le
 \frac{2i\pi}{k}}S^-_M(e^{\sqrt{-1}\theta})\right)\nn\\
 &&=m(i(\gamma^2)+S^+_M(1))+\sum_{\theta\in(0,\pi)}\left(\sum_{\frac{k\theta}{2\pi}< i\le m}
 S^+_M(e^{\sqrt{-1}\theta})-
 \sum_{\frac{k\theta}{2\pi}\le i\le
 m}S^-_M(e^{\sqrt{-1}\theta})\right)\nn\eea
  \bea&&=m(i(\gamma^2)+S^+_M(1))+\sum_{\theta\in(0,\pi)}\left(\left(m
 -\left[\frac{k\theta}{2\pi}\right]\right)S^+_M(e^{\sqrt{-1}\theta})-\left[m+1
 -\frac{k\theta}{2\pi}\right]S^-_M(e^{\sqrt{-1}\theta})\right)\nn\\
&&=m(i(\gamma^2)+S^+_M(1))\nn\\
&&\quad +\sum_{\theta\in(0,\pi)}\left(\left(m
 -\left[\frac{k\theta}{2\pi}\right]\right)S^-_M(e^{\sqrt{-1}(2\pi-\theta)})
-\left(m+1
 -E\left(\frac{k\theta}{2\pi}\right)\right)S^-_M(e^{\sqrt{-1}\theta})\right)\nn\\
&&=m(i(\gamma^2)+S^+_M(1))+\sum_{\theta\in(\pi,2\pi)}\left(m
 -\left[\frac{k(2\pi-\theta)}{2\pi}\right]\right)S^-_M(e^{\sqrt{-1}\theta})\nn\\
&&\quad-\sum_{\theta\in(0,\pi)}\left(m+1
 -E\left(\frac{k\theta}{2\pi}\right)\right)S^-_M(e^{\sqrt{-1}\theta})\nn\\
 &&= m(i(\gamma^2)+S^+_M(1))+\sum_{\theta\in(0,\pi)\cup
 (\pi,2\pi)}\left(-(m+1)+E\left(\frac{k\theta}{2\pi}\right)\right)S^-_M(e^{\sqrt{-1}\theta})\nn\\
 &&=m(i(\gamma^2)+S^+_M(1))-(m+1)C(M)+\sum_{\theta\in
 (0,2\pi)}E\left(\frac{k\theta}{2\pi}\right)S^-_M(e^{\sqrt{-1}\theta})\nn\\
 &&=m(i(\gamma^2)+S^+_M(1)-C(M))
+\sum_{\theta\in(0,2\pi)}E\left(\frac{k\theta}{2\pi}\right)S_M^-(e^{\sqrt{-1}\theta})-C(M),\nn
 \eea
where in the fourth equality and sixth equality we have used the
facts that
$$S_M^+(e^{\sqrt{-1}\theta})=S_M^-(e^{\sqrt{-1}(2\pi-\theta)}),\;\;
$$
$k=2m+1$ and $E(a)+[b]=a+b$ if $a,\; b\in \R$ and $a+b\in \Z$,
especially $E(-a)+[a]=0$ for any $a\in \R$. By using Theorem 4.1 and
$m=\frac {k-1}{2}$ we get (\ref{1.21}).
 Similarly we obtain (\ref{1.23}).
\hfill\hb

 \noindent{\bf Corollary 5.1.} {\it For mean $L_0$-index, there holds
$$\hat{i}_{L_0}(\gamma)=\frac{1}{2}\hat{i}(\ga^2)
=\frac 12
(i(\gamma^2)+S^+_M(1)-C(M))+\sum_{\theta\in(0,2\pi)}\frac{\theta}{2\pi}S^-_M(e^{\sqrt{-1}\theta}).
$$}
\noindent{\bf Proof.} The above equality follows from Theorem 5.1
and the definition of the mean $L_0$-index
$$\hat {i}_{L_0}(\gamma)=\lim_{k\to \infty}\frac{i_{L_0}(\gamma^k)}{k}.$$
\hfill\hb

\noindent{\bf Proposition 5.1.}(Theorem 4.3 in \cite{LZ})  {\it Let
$\ga_j\in \mathcal {P}_{\tau_j}(2n)$ for $j=1,\cdots,q$ be a finite
collection of
 Symplectic paths. Extend $\ga_j$ to $[0,+\infty)$
 by $\ga_j(t+\tau_j)=\ga_j(t)\ga_j(\tau_j)$ and let
 $M_j=\ga(\tau_j)$,  for $j=1,\cdots,q$ and $t>0$.
 Suppose \bea \hat{i}(\ga_j)>0, \quad
          j=1,\cdots,q.\nn\eea
  Then there exist infinitely many $(R, m_1, m_2,\cdots,m_q)\in \N^{q+1}$ such that

  (i) $\nu(\ga_j, 2m_j\pm 1)=\nu(\ga_j)$,

  (ii) $i(\ga_j, 2m_j-1)+\nu(\ga_j,
  2m_j-1)=2R-(i(\ga_j)+2S_{M_j}^+(1)-\nu(\ga_j))$,

  (iii)$i(\ga_j,2m_j+1)=2R+i(\ga_j)$,

\noindent where we have set $i(\ga_j, n_j)=i(\ga_j, [0,
n_j\tau_j])$, $\nu(\ga_j, n_j)=\nu(\ga_j, [0, n_j\tau_j])$ for
$n_j\in\N$.}

\noindent{\bf Proof of Theorem 1.5.} We divide our proof in three
steps.

\noindent{\it Step 1.} Application of Proposition 5.1.

 By (\ref{6.11}) and (\ref{6.12}), we have
    \be \hat{i}(\ga_j^2)=2\hat{i}_{L_0}(\ga_j)>0.\lb{6.13}\ee
So we have
   \be \hat{i}(\ga_j^2)>0, \quad
          j=1,\cdots,q,\lb{6.14}\ee
where $\ga_j^2$ is the 2-times iteration of $\ga_j$ defined by
(\ref{4.4}).
  Hence the symplectic paths $\ga_j^2, j=1,2,\cdots,q$ satisfy the condition in Theorem
  6.1,
so there exist infinitely $(R, m_1, m_2,\cdots,m_q)\in \N^{q+1}$
such that
\bea \nu(\ga_j^2, 2m_j\pm 1)&=&\nu(\ga_j^2),\lb{6.15}\\
 i(\ga_j^2, 2m_j-1)+\nu(\ga_j^2,
 2 m_j-1)
  &=&2R-(i(\ga_j^2)+2S_{M_j}^+(1)-\nu(\ga_j^2)),\lb{6.16}\\
 i(\ga_j^2,2m_j+1)&=&2R+i(\ga_j^2).\lb{6.17}\eea

\noindent{\it Step 2.} Verification of (i).

By  Theorems 4.1 and 4.2, we have
 \bea
 \nu_{L_0}(\ga_j,2m_j\pm 1)=\nu_{L_0}(\ga_j)+\frac{\nu(\ga_j^2,2m_j\pm 1)-\nu(\ga^2_j)}{2},\lb{6.20}\\
 \nu_{L_1}(\ga_j,2m_j\pm 1)=\nu_{L_1}(\ga_j)+\frac{\nu(\ga_j^2,2m_j\pm 1)-\nu(\ga^2_j)}{2}. \lb{6.20'}\eea
Hence (i) follows from (\ref{6.15}) and (\ref{6.20}).

\noindent{\it Step 3.} Verifications of (ii) and (iii).

 By Theorems 4.1 and 4.2, we have
 \bea
&& i_{L_0}(\ga^m)-i_{L_1}(\ga^m)=i_{L_0}(\ga)-i_{L_1}(\ga),
\quad \forall m\in 2\N-1,\lb{6.21}\\
&&i_{L_0}(\ga^m)-i_{L_1}(\ga^m)=i_{L_0}(\ga^2)-i_{L_1}(\ga^2), \quad
\forall m\in 2\N.\lb{6.22}\eea
 By (\ref{6.10}), (\ref{6.9})  and (\ref{6.21}) we have \be
2i_{L_0}(\ga_j,2m_j\pm 1)=i(\ga_j^2,2m_j\pm
1)-n+i_{L_0}(\ga_j)-i_{L_1}(\ga_j).\lb{6.23}\ee
 By (\ref{6.15}),
(\ref{6.16}) and (\ref{6.23}) we have
 \be
 2i_{L_0}(\ga_j,,2m_j-1)=2R-(i(\ga_j^2)-2S_{M_j}^+(1)+n-i_{L_0}(\ga_j)+i_{L_1}(\ga_j)).\lb{6.24}\ee
So by (\ref{6.10}) we have
 \be
i_{L_0}(\ga_j,2m_j-1)=R-(i_{L_1}(\ga_j)+n+S_{M_j}^+(1)).\lb{6.25}\ee
Together with (i), this yields (ii).

By (\ref{6.17}) and (\ref{6.23}) we have
     \be
     2i_{L_0}(\ga_j,2m_j+1)=2R+i(\ga_j^2)-n+i_{L_0}(\ga_j)-i_{L_1}(\ga_j).\lb{6.26}\ee
By (\ref{6.10}) and (\ref{6.26}) we have \be
i_{L_0}(\ga_j,2m_j+1)=R+i_{L_0}(\ga_j).\lb{6.27}\ee Hence (iii)
holds and the proof of Theorem 1.5 is complete. \hfill\hb

\noindent{\bf Remark 5.1.} From (\ref{6.12}) and (iii) of
 Theorem 1.5, it is easy to see that for any $\mathcal {R}>0$, among
the infinitely many vectors $(R,m_1,m_2,\cdots,m_q)\in \N^{q+1}$ in
Theorem 1.5, there exists one vector such that its first component
$R$ satisfies $R>\mathcal {R}$.

\setcounter{equation}{0}
\section{Variational set up}

In this section, we briefly recall the variational set up and some
corresponding results proved in \cite{LZZ}. Based on these results
we obtain an injection map in Lemma 6.3 bellow which is  basic in
the proofs of Theorems 1.1 and 1.2.

   For $\Sg\in \mathcal{H}_b^{s,c}(2n)$, let $j_\Sg: \Sg \rightarrow[0,+\infty)$ be the
gauge function of $\Sg$ defined by
        \be
        j_{\Sg}(0)=0,\quad {\rm and} \quad  j_\Sg(x)=\inf\{\lambda >0\mid
        \frac{x}{\lambda}\in C\}, \quad \forall x \in
        \R^{2n}\setminus\{0\},\lb{7.1}
        \ee
where $C$ is the domain enclosed by $\Sg$.

 Define
        \be H_\alpha(x)=(j_\Sg(x))^\alpha,\;\alpha>1,\quad
        H_\Sg(x)=H_2(x),\; \forall x \in
        \R^{2n}.\lb{7.2}
        \ee
Then $H_{\Sigma} \in C^2 (\R^{2n}\backslash \{0\},\R)\cap
    C^{1,1}(\R^{2n},\R)$. Its Fenchel conjugate (cf.\cite{EH},\cite{Ek}) is the function
    $H_\Sg^*$ defined by
        \be
          H_\Sg^*(y)=\max\{(x\cdot y
          -H_\Sg(x))|\, x\in \R^{2n}\}.\lb{7.3}\ee

We consider the following fixed energy problem
\bea
\dot{x}(t) &=& JH_\Sg'(x(t)), \lb{7.4}\\ H_\Sg(x(t)) &=& 1,   \lb{7.5}\\
x(-t) &=& Nx(t),  \lb{7.6}\\ x(\tau+t) &=& x(t),\quad \forall\,
t\in\R. \lb{7.7} \eea \
 Denote by
$\mathcal{J}_b(\Sg,2)\;(\mathcal{J}_b(\Sg,\alpha)$ for $\alpha=2$ in
(\ref{7.2})) the set of all solutions $(\tau,x)$ of problem
(\ref{7.4})-(\ref{7.7}) and by $\tilde{\mathcal{J}}_b(\Sg,2)$ the
set of all geometrically distinct solutions of
(\ref{7.4})-(\ref{7.7}). By Remark 1.2 or discussion in \cite{LZZ},
elements in $\mathcal{J}_b(\Sg)$ and $\mathcal{J}_b(\Sg,2)$ are one
to one correspondent. So we have
$^\#\td{\mathcal{J}}_b(\Sg)$=$^\#\td{\mathcal{J}}_b(\Sg,2)$.

For $S^1=\R / \Z$, as in \cite{LZZ} we define the Hilbert space $E$
by
    \bea E = \left\{ x\in W^{1,2}(S^1,\R^{2n})\left| x(-t)=Nx(t),\quad {\rm for\; all} \;
    t\in \R \;\;{\rm and} \;\; \int_0^1x(t)dt=0 \right.\right\}. \lb{7.8} \eea
The inner product on $E$ is given by
   \be
    (x,y)=\int_0^1 \< \dot{x}(t), \dot{y}(t) \> dt.\lb{7.9}
    \ee

   The $C^{1,1}$ Hilbert manifold $M_\Sg \subset E$ associated to $\Sg$ is
defined by
  \be M_\Sg=\left\{ x\in E \left| \int_0^1H^*_\Sg(-J\dot{x}(t))dt=1\; {\rm and} \;
  \int_0^1\< J\dot{x}(t), x(t)\>dt <0\right.\right\}. \lb{7.10}\ee

  Let $\Z_2=\{-id, id\}$ be the usual $\Z_2$ group. We define the $\Z_2$-action on $E$ by
        $$-id(x)=-x,\quad id(x)=x, \qquad \forall x\in E.$$
Since $H^*_\Sg$ is even, $M_\Sg$ is symmetric to 0, i.e., $\Z_2$
invariant. $M_\Sg$ is a paracompact $\Z_2$-space. We
  define
     \be
       \Phi(x)=\frac{1}{2}\int_0^1\< J\dot{x}(t), x(t)\>dt, \lb{7.11}
       \ee
  then $\Phi$ is a $\Z_2$ invariant function  and $\Phi\in C^\infty (E,\R)$. We denote by $\Phi_\Sg$
  the restriction of $\Phi$ to $M_\Sg$, we remind that   $\Phi$ and
  $\Phi_{\Sg}$ here
  are the functionals $A$ and $A_{\Sg}$ in \cite{LZZ} respectively.

Suppose $z\in M_\Sg$ is a critical point of $\Phi_\Sg$. By Lemma 7.1
of \cite{LZZ} there is a $c_1(z)\in 0\times \R^n$ such that
$x(z)(t)=(|\Phi_\Sg(z)|^{-1}(z(|\Phi_\Sg(z)|t)+c_1(z))$ is a
$\tau$-periodic solution of the fixed energy problem
(\ref{1.11})-(\ref{1.12}), i.e., $(\tau ,x)\in \mathcal{J}_b(\Sg,2)$
with $\tau=|\Phi_\Sg(z)|^{-1}$.

Following  the ideas of Ekeland and Hofer in \cite{EH},  Long, Zhu
and the second author of this paper in  \cite{LZZ} proved the
following result(see Corollary 7.10 of \cite{LZZ}).

\noindent{\bf Lemma 6.1.} {\it If
$^\#\td{\mathcal{J}}_b(\Sg)<+\infty$, then for each $k\in \N$, there
exists a critical points $z_k\in M_\Sg$ of $\Phi_\Sg$ such that the
sequence $\{\Phi_{\Sg}(z_k)\}$ increases strictly to zero as $k$
goes to $+\infty$ and there holds
          $$m^{-}(z_k)\le k-1\le m^-(z_k)+m^0(z_k),$$
where $m^{-}(z_k)$ and $m^0(z_k)$ are Morse index and nullity of the
formal Hessian $Q_{z_k}$ of $\Phi_\Sg$ at $z$ defined by (7.36) of
\cite{LZZ} as follows:
  \be Q_{z_k}(h)=\frac {1}{2}\int_0^1\langle
  J\dot{h}(t),h(t)\rangle dt-\frac {1}{2}\Phi(z_k)\int_0^1\langle
  (H^{*}_\Sg)''(-J\dot {z}_k(t))J\dot{h}(t),J\dot{h}(t)\rangle dt,\;\;h\in
  T_{z_k}M_\Sg.\lb{7.14}\ee}

We remind that $L_0=\{0\}\times\R^n$ and $L_1=\R^n\times\{0\}\subset
\R^{2n}$. The following two maslov-type indices are defined in
\cite{LZZ}.

\noindent{\bf Definition 6.1.} {\it For
$M=\left(\begin{array}{cc}A&B\\C&D\end{array}\right)\in \Sp(2n)$, we
define
  \be \nu_1(M)=\dim \ker B,\quad
{\rm and}\quad
      \nu_2(M)=\dim \ker
    C.\lb{6.1}\ee
For $\Psi\in C([a,b],\Sp(2n))$, we define
  \be \nu_1(\Psi)=\nu_1 (\Psi(b)),\quad
\quad
      \nu_2(\Psi)=\nu_2 (\Psi(b))\lb{6.2}\ee
and
   \be \mu_1(\Psi,[a,b])=i_{{CLM}_{\R^{2n}}}(L_0, \Psi L_0,
   [a,b]),\quad \mu_2(\Psi,[a,b])=i_{{CLM}_{\R^{2n}}}(L_1, \Psi L_1,
   [a,b]),\lb{6.3}\ee
where the Maslov index $i_{{CLM}_{\R^{2n}}}$ for Lagrangian subspace
paths is defined in \cite{CLM}. We will omit the interval $[a,b]$ in
the index notations when there is no confusion.}

 By Proposition
C of \cite{LZZ}, we have
 \bea \mu_1(\ga)+\mu_2(\ga)=i(\ga^2)+n,\quad
\nu_1(\ga)+\nu_2(\ga)=\nu(\ga^2),\lb{6.10}\eea where $\ga^2$ is the
2-times iteration of $\ga$ defined by (\ref{4.4}).

For convenience in the further proofs of Theorems 1.1 and 1.2 in
this paper, we firstly give a relationship between the Maslov-type
indices $\mu_1$, $\mu_2$ and $i_{L_0}$, $i_{L_1}$.

 \noindent{\bf Proposition
6.1.} {\it For any $\ga\in \mathcal{P}_\tau(2n)$, there hold \bea
&&\nu_1(\ga)=\nu_{L_0}(\ga),\quad
       \nu_2(\ga)=\nu_{L_1}(\ga),\lb{6.8}\\
&&\mu_1(\ga)=i_{L_0}(\ga)+n, \quad
 \mu_2(\ga)=i_{L_1}(\ga)+n.\lb{6.9}\eea
}

From (\ref{4.24}) and (\ref{6.10})-(\ref{6.9}), we have \bea
\hat{\mu}_1(\ga)=\hat{\mu}_2(\ga)=\hat{i}_{L_0}(\ga)=\hat{i}_{L_1}(\ga)=\frac{1}{2}\hat{i}(\ga^2),\lb{6.11}\eea
where $\hat{\mu}_j(\ga)$ is the $\mu_j$-mean index for $j=1,2$
defined in \cite{LZZ}.

\noindent{\bf Proof.} (\ref{6.8}) follows from the definitions of
$\nu_{L_0}$ and $\nu_{L_1}$ in Definitions 2.1 and 2,4 and the
definitions of $\nu_1$ and $\nu_2$ in Definitions 6.1.

 (\ref{6.9}) follows from (\ref{c1}) and Theorem 2.4 of
 \cite{Zhang2}. We note that for $x,y\in W_1$, there hold
 $$(Ax,y)=2(A^1x,y),\;\;(Bx,y)=2(B^1x,y),$$
 where $W_1$, $A,\;B$ were defined in \cite{Zhang2} before Theorem
 2.4.
\hfill\hb

 By Proposition 5.1,
 Lemma 8.3 of \cite{LZZ} and Lemma 6.1, we have the following
result which is also basic in the proof of Theorems 1.1 and 1.2.

\noindent{\bf Lemma 6.2.} {\it If
$^\#\tilde{\mathcal{J}}_b(\Sg)<+\infty$, there is an sequence
$\{c_k\}_{k\in \N}$, such that \bea
-\infty<c_1<c_2<\cdots<c_k<c_{k+1}<\cdots<0,\lb{7.15}\\
c_k\rightarrow 0\quad {\rm as}\;k\rightarrow +\infty.\lb{7.16}\eea
For any $k\in \N$, there exists a brake orbit $(\tau, x)\in
\mathcal{J}_b(\Sg,2)$ with $\tau$ being the minimal period of $x$
and $m\in \N$ satisfying $m\tau=(-c_k)^{-1}$ such that
 for
     \be
     z(x)(t)=(m\tau)^{-1}x(m\tau t)-\frac{1}{(m\tau)^2}\int_0^{m\tau}
     x(s)ds, \quad t\in S^1 , \lb{7.17}\ee
$z(x)\in M_\Sg$ is a critical point of $\Phi_\Sg$ with
$\Phi_\Sg(z(x))=c_k$ and \be i_{L_0}(x,m)\le k-1\le
i_{L_0}(x,m)+\nu_{L_0}(x,m)-1,\lb{z}\lb{7.18}\ee
 where we denote by $(i_{L_0}(x,m),\nu_{L_0}(x,m))=(i_{L_0}(\ga_x,m),\nu_{L_0}(\ga_x,m))$ and $\ga_x$ the associated
 symplectic path of $(\tau,x)$.  }

 \noindent{\bf
Definition 6.2.} {\it We call $(\tau,x)\in\mathcal{J}_b(\Sg,2)$ with
minimal period $\tau$ {\it infinitely variational visible} if there
are infinitely many $m's\in \N$ such that $(\tau,x)$ and $m$ satisfy
conclusions in Lemma 6.2. We denote by $\mathcal
{V}_{\infty,b}(\Sg,2)$ the subset of $\td{\mathcal{J}}_b(\Sg,2)$
consisting of $[(\tau,x)]$ in which there is an infinitely
variational visible representative.}

We have the following injective map lemma about the $L_0$-index.

\noindent {\bf Lemma 6.3.} {\it Suppose
$^\#\tilde{\mathcal{J}}_b(\Sg)<+\infty$. Then there exist an integer
$K\ge 0$ and an injection map $\phi: \N+K\mapsto
\mathcal{V}_{\infty,b}(\Sg,2)\times \N$ such that

(i) For any $k\in \N+K$, $[(\tau,x)]\in
\mathcal{V}_{\infty,b}(\Sg,2)$ and $m\in \N$ satisfying
$\phi(k)=([(\tau \;,x)],m)$, there holds
              $$i_{L_0}(x,m)\le k-1\le i_{L_0}(x,m)+\nu_{L_0}(x,m)-1,$$
where $x$ has minimal period $\tau$.

(ii) For any $k_j\in \N+K$, $k_1<k_2$, $(\tau_j,x_j)\in \mathcal
{J}_b(\Sg,2)$ satisfying $\phi(k_j)=([(\tau_j \;,x_j)],m_j)$ with
$j=1,2$ and $[(\tau_1 \;,x_1)]=[(\tau_2 \;,x_2)]$, there holds
       $$m_1<m_2.$$}
\noindent{\bf Proof.} Since $^\#\tilde{\mathcal{J}}_b(\Sg)<+\infty$,
there is an integer $K\ge 0$ such that all critical values $c_{k+K}$
with $k\in \N$ come from iterations of elements in
$\mathcal{V}_{\infty,b}(\Sg,2)$. Together with Lemma 6.2, for each
$k\in \N$, there is a $(\tau,x)\in \mathcal{J}_{b}(\Sg,2)$ with
minimal period $\tau$ and $m\in \N$ such that (\ref{7.17}) and
(\ref{7.18}) hold for $k+K$ instead of $k$. So we define a map
$\phi:\N+K\mapsto \mathcal {V}_{\infty,b}(\Sg,2)\times \N$ by
$\phi(k+K)=([(\tau,x)],m)$.

For any $k_1<k_2\in \N$, if $\phi(k_j)=([\tau_j,x_j)],m_j)$ for
$j=1,2$. Write $[(\tau_1,x_1)]=[(\tau_2,x_2)]=[(\tau,x)]$ with
$\tau$ being the minimal period of $x$, then by Lemma 6.2 we have
        \be m_j\tau=(-c_{k_j+K})^{-1},\quad j=1,2.\lb{7.19}\ee
 Since $k_1<k_2$ and $c_k$ increases strictly to 0 as $k\rightarrow
 +\infty$, we have
          \be m_1<m_2.\lb{7.20}\ee
 So the map $\phi$ is injective, also (ii) is proved.
 The proof of this Lemma 6.3 is complete.
 \hfill\hb
\setcounter{equation}{0}
\section{Proof of Theorem 1.1} 

 We first prove Lemma 1.1.

\noindent{\bf Proof of Lemma 1.1.} We set
$\ga(\frac{\tau}{2})=\left(\begin{array}{cc}A&B\\C&D\end{array}\right)$
in square block form.
 Since $(\tau ,x)\in
\mathcal{J}_b(\Sg,2)$, we have
             \be \dot{x}(t)=JH'_\Sg(x(t)),\quad t\in \R.\lb{8.4}\ee
  By the definition of $H_\Sg$ in (\ref{7.2}),
    $H_\Sg$ is 2-homogeneous and $H'_\Sg$ is 1-homogeneous . So we have
    \be \dot{x}(t)=JH_\Sg''(x(t))x(t),\quad t\in \R.\lb{8.5}\ee
Differentiating (\ref{8.4}) we obtain
     \be \ddot{x}(t)=JH_\Sg''(x(t))\dot{x}(t),\quad t\in \R.\lb{8.6}\ee
Since $\ga$ is the associated symplectic path of $(\tau,x)$,
   $\ga(t)$ is the solution of the problem
   \bea
   \dot{\ga}(t) &=& JH_\Sg''(x(t))\ga(t), \lb{8.7}\\
   \ga(0) &=& I_{2n}.  \lb{8.8} \eea
So we have
   \be x(t)=\ga(t)x(0),\quad \dot{x}(t)=\ga(t)\dot{x}(0), \qquad t\in
   \R.\lb{8.9}\ee
Denote by $x(t)=(p(t),q(t))\in \R^n\times\R^n$. Since
 \be x(-t)=Nx(t),\quad x(t+\tau)=x(t),\qquad t\in \R, \lb{8.10}\ee
 we have
 \bea p(0)=0=p(\frac{\tau}{2}), \;q(0)\neq 0,\lb{8.11}\\
      \dot{p}(0)\neq 0,\; \dot{q}(0)=0=\dot{q}(\frac{\tau}{2}).\lb{8.12}\eea
Since $(\tau,x)$ is symmetric, by (\ref{8.9}) we have
   \bea
   \left(\begin{array}{c}0\\-q(0)\end{array}\right)&=&\left(\begin{array}{c}0\\q(\frac{\tau}{2})\end{array}\right)=
   \left(\begin{array}{c}p(\frac{\tau}{2})\\q(\frac{\tau}{2})\end{array}\right)=
   \left(\begin{array}{cc}A&B\\C&D\end{array}\right)\left(\begin{array}{c}p(0)\\q(0)\end{array}\right)\nn\\
   &=&\left(\begin{array}{cc}A&B\\C&D\end{array}\right)\left(\begin{array}{c}0\\q(0)\end{array}\right)
   =\left(\begin{array}{c}Bq(0)\\Dq(0)\end{array}\right),\\\nn\\
   \left(\begin{array}{c}-\dot{p}(0)\\0\end{array}\right)&=&\left(\begin{array}{c}\dot{p}(\frac{\tau}{2})\\0\end{array}\right)=
   \left(\begin{array}{c}\dot{p}(\frac{\tau}{2})\\\dot{q}(\frac{\tau}{2})\end{array}\right)=
   \left(\begin{array}{cc}A&B\\C&D\end{array}\right)\left(\begin{array}{c}\dot{p}(0)\\
   \dot{q}(0)\end{array}\right)\nn\\
   &=&\left(\begin{array}{cc}A&B\\C&D\end{array}\right)\left(\begin{array}{c}\dot{p}(0)\\0\end{array}\right)
   =\left(\begin{array}{c}A\dot{p}(0)\\C\dot{p}(0)\end{array}\right). \eea
So we have
    \bea &&B q(0)=0,\quad C\dot{p}(0)=0,\lb{8.15}\\
         &&D q(0)=-q(0),\quad A\dot{p}(0)=-\dot{p}(0). \lb{8.16}\eea
Since
     \be \langle Jx(0),\dot{x}(0)\rangle= \langle
     Jx(0),JH'_\Sg(x(0))\rangle=\langle
     x(0),H'_\Sg(x(0))\rangle=2H_\Sg(x(0))=2,\lb{8.17}\ee
where we have used the fact that $(\tau,x)\in \mathcal{J}_b(\Sg,2)$
and $H_\Sg$ is 2-homogeneous, we have
     \be \langle q(0), \dot{p}(0)\rangle=-\langle
     Jx(0),\dot{x}(0)\rangle=-2.\lb{8.18}\ee
    Denote by $\xi=-\frac{1}{\sqrt{2}}\dot{p}(0)$ and $\eta=\frac{1}{\sqrt{2}}q(0)$.
    We have
    \be \xi^T\eta=1, \lb{8.19}\ee
and
\bea &&B\eta=0,\quad C\xi=0,\lb{8.20}\\
         &&D\eta=-\eta,\quad A\xi=-\xi, \lb{8.21}\eea
  where we denote by $\xi^T$ the transpose of $\xi$.

\noindent{\it Claim.} There exist two $n\times (n-1)$ matrices $F$
    and $G$  such that $\det(\xi F)>0$ and the matrix $\left(\begin{array}{cc} (\xi
    F)
    &0\\0&(\eta G)\end{array}\right)\in \Sp(2n)$, where $(\xi
    F)$ and $(\eta G)$ are $n\times n$ matrices whose first columns are
    $\xi$ and $\eta$, and the other $n-1$ columns are the matrices $F$ and $G$ respectively.

\noindent{\it Proof of the claim.} We divide the proof into two
cases.

\noindent {\it Case 1.}  $\xi=\lm \eta$ for some $\lm\in
\R\setminus\{0\}$. Denote by $\span \{e_2,e_3,\cdots,e_n\}$ the
orthogonal complement of $\span\{\xi\}$ in $\R^n$ in the standard
inner product sense, where $e_2, e_3,\cdots,e_n$ are unit and mutual
orthogonal. Define the $n\times(n-1)$ matrix
$\td{F}=(e_2\;e_3\;\cdots\;e_n)$ whose columns are
$e_2,e_3,\cdots,e_n$. If $\det (\xi \tilde{F})>0$, we define
$F=G=(e_2\;e_3\;\cdots\;e_n)$. Otherwise we define
$F=G=\left((-e_2)\;e_3\;e_4\;\cdots\;e_n\right)$. By direct
computation we always have $\det(\xi F)>0$ and the matrix
$\left(\begin{array}{cc}(\xi F)
    &0\\0&(\eta G)\end{array}\right)\in \Sp(2n)$.

\noindent {\it Case 2.}   $\xi\neq \lm \eta$ for all $\lm\in
\R\setminus\{0\}$, i.e., $\dim \span\{\xi,\eta\}=2$. Denote by
$\span \{e_3,\cdots,e_n\}$ the orthogonal complement of
$\span\{\xi,\eta\}$ in $\R^n$ in the standard inner product sense,
where $ e_3,\cdots,e_n$ are unit and mutual orthogonal. Denote by $
\span\{\xi,\eta\}=\span\{e_1,e_2\}$ where  $e_1$ and $e_2$ are unit
and orthogonal and $\lm e_1=\xi$ for some $\lm\in \R$. Since
$\xi^T\eta=1$ we have $\eta=\lm^{-1}e_1+re_2$ for some $r\in
\R\setminus\{0\}$. Then we define the matrix $\td{F}=((\lm
e_1-r^{-1}e_2)\; e_3\;.\;.\;.\;e_n)$ whose columns are $ \lm
e_1-r^{-1}e_2,\;e_3,\cdots,e_n$. If $\det(\xi\,\td{F})>0$, we define
$F=((\lm e_1-r^{-1}e_2)\; e_3\;e_4\;.\;.\;.\;e_n)$ and $G=(
(-re_2)\; e_3\;e_4\;.\;.\;.\;e_n)$. Otherwise we define $F=((\lm
e_1-r^{-1}e_2)\; e_3\;.\;.\;.\;(-e_n))$ and $G=(-re_2\;
e_3\;e_4\;.\;.\;.\;(-e_n))$. By direct computation we always have
$\det(\xi F)>0$ and the matrix $\left(\begin{array}{cc}(\xi F)
    &0\\0&(\eta G)\end{array}\right)\in \Sp(2n)$.
By the discussion in cases 1 and 2, the claim  is proved.

By this claim, there exist two $n\times (n-1)$ matrices $F$
    and $G$  such that $\det(\xi F)>0$ and the matrix $\left(\begin{array}{cc}(\xi
    F)
    &0\\0&(\eta G)\end{array}\right)\in \Sp(2n).$ So we have
    \be (\eta G)=((\xi F)^T)^{-1}.\lb{8.22}\ee
    Applying (\ref{8.20})-(\ref{8.22}), by direct
    computation we have
    \bea && \left(\begin{array}{cc}(\eta G)^T
    &0\\0&(\xi F)^T\end{array}\right)\left(\begin{array}{cc}A
    &B\\C&D\end{array}\right)\left(\begin{array}{cc}(\xi F)
    &0\\0&(\eta G)\end{array}\right)\nn\\
    &=& \left(\begin{array}{cccc}-1
    &\eta^TAF&0&\eta^TBG\\0&G^TAF&0&G^TBG\\0&\xi^TCF&-1&\xi^TDG\\
    0&F^TCF&0&F^TDG\end{array}\right).\lb{8.23}\eea
Since the above matrix is still a symplectic matrix, by Lemma 1.1.2
of \cite{Long1}, we have that both $\left(\begin{array}{cc}-1
    &0\\(\eta^TAF)^T&(AF)^TG\end{array}\right)\left(\begin{array}{cc}0
    &\xi^TCF\\0&F^TCF\end{array}\right)$ and
     $ \left(\begin{array}{cc}0
    &0\\(\eta^TBG)^T&G^TB^TG
    \end{array}\right)\left(\begin{array}{cc}-1&\xi^TDG
    \\ 0&F^TDG\end{array}\right)$
     are  symmetric and
     \bea \left(\begin{array}{cc}-1
    &0\\(\eta^TAF)^T&(AF)^TG\end{array}\right)\left(\begin{array}{cc}-1
    &\xi^TDG\\0&F^TDG\end{array}\right)
    -\left(\begin{array}{cc}0
    &0\\(\xi^T(CF))^T&(CF)^TF\end{array}\right)\left(\begin{array}{cc}0
    &\eta^TBG\\0&G^TBG\end{array}\right)=I_n.\nn\eea
So by the above three facts and direct computation we have \bea
\eta^TAF=0, \quad \eta^TBG=0,\quad
      \xi^TCF=0,\quad \xi^TDG=0.\lb{8.25}\eea
Set $\td{M}=\left(\begin{array}{cc}G^TAF
    &G^TBG\\F^TCF&F^TDG\end{array}\right)$. By
    (\ref{8.23}) and (\ref{8.25}), there hold
 $\td{M}\in \Sp(2n-2)$ and

 \bea  \left(\begin{array}{cc}(\eta G)^T
    &0\\0&(\xi F)^T\end{array}\right)\left(\begin{array}{cc}A
    &B\\C&D\end{array}\right)\left(\begin{array}{cc}(\xi F)
    &0\\0&(\eta G)\end{array}\right)
   = (-I_2)\diamond \tilde{M}.\lb{8.26}\eea
Since $\det(\xi F)>0$, there is a continuous matrix path $\psi(s)$
for $s\in [0,1]$ joints $(\xi F)$ and $I_n$ such that $\psi(0)=I_n$
and $\psi(1)=(\xi F)$ and $\det(\psi(s))>0$ for all $s\in [0,1]$.
For $s\in [0,1]$, we define \be
\Psi(s)=\left(\begin{array}{cc}\psi(s)^{-1}&0\\0&\psi(s)^T\end{array}\right)
          \left(\begin{array}{cc}A&B\\C&D\end{array}\right)
          \left(\begin{array}{cc}\psi(s)&0\\0&(\psi(s)^T)^{-1}\end{array}\right).\lb{8.27}\ee
Then by (\ref{8.22}) and (\ref{8.26}),  $\Psi$ satisfies the
conclusions in Lemma 1.1 and the proof is complete. \hfill\hb

 In order to prove Theorem 1.1, we need the following three results.

 \noindent{\bf Lemma 7.1.} {\it For any symmetric
$(\tau ,x)\in \mathcal{J}_b(\Sg,2)$, denote by $\ga$ the symplectic
path associated to $(\tau,x)$. We have
 \be
 \left|\left(i_{L_0}(\ga)+\nu_{L_0}(\ga)\right)-
 \left(i_{L_1}(\ga)+\nu_{L_1}(\ga)\right)\right|\le
 n-1.\lb{8.28}\ee}
\noindent{\bf Proof.}
 By Lemma 1.1 there exist a symplectic path
$\ga^*\in \mathcal{P}_\frac{\tau}{2}(2n)$ and $\td{M}\in \Sp(2n-2)$
such that \be \ga\; \sim_{L_j}\; {\ga^*} \qquad {\rm for}\quad
j=0,\;1,\lb{8.29}\ee
 \be {\ga^*}(\frac{\tau}{2})=(-I_2)\diamond
\td{M}.\lb{8.30}\ee
 So by Theorem 2.1, we have
 \bea
&&\left|\left(i_{L_0}(\ga)+\nu_{L_0}(\ga)\right)-
 \left(i_{L_1}(\ga)+\nu_{L_1}(\ga)\right)\right|\nn\\
&=&\left|\left(i_{L_0}(\ga^*)+\nu_{L_0}(\ga^*)\right)-
 \left(i_{L_1}(\ga^*)+\nu_{L_1}(\ga^*)\right)\right|.
\lb{8.31}\eea
 We choose a special symplectic
path $\td{\ga}=\ga_1\diamond \ga_2\in
\mathcal{P}_{\frac{\tau}{2}}(2n)$, where $\ga_1\in
\mathcal{P}_{\frac{\tau}{2}}(2)$, $\ga_1({\frac{\tau}{2}})=-I_2$ and
$\ga_2\in \mathcal{P}_{\frac{\tau}{2}}(2n-2)$,
$\ga_2({\frac{\tau}{2}})=\td{M}$.

 By Theorems 2.2 and 2.3, we have
\bea &&\left|\left(i_{L_0}(\ga^*)+\nu_{L_0}(\ga^*)\right)-
 \left(i_{L_1}(\ga^*)+\nu_{L_1}(\ga^*)\right)\right|\nn\\
&=&\left|\left(i_{L_0}(\td{\ga})+\nu_{L_0}(\td{\ga})\right)-
 \left(i_{L_1}(\td{\ga})+\nu_{L_1}(\td{\ga})\right)\right|\nn\\
&=&|\left(i_{L_0}(\ga_1)+\nu_{L_0}(\ga_1)\right)-
 \left(i_{L_1}(\ga_1)+\nu_{L_1}(\ga_1)\right)\nn\\
 &&\;+\left(i_{L_0}(\ga_2)+\nu_{L_0}(\ga_2)\right)-
 \left(i_{L_1}(\ga_2)+\nu_{L_1}(\ga_2)\right)|.\lb{8.32}\eea
Since $-I_2\in O(2)\cap \Sp(2)$, by Theorem 2.3 again we have \bea
&&\left(i_{L_0}(\ga_1)+\nu_{L_0}(\ga_1)\right)-
 \left(i_{L_1}(\ga_1)+\nu_{L_1}(\ga_1)\right)=0,\lb{8.33}\\
&&|\left(i_{L_0}(\ga_2)+\nu_{L_0}(\ga_2)\right)-
 \left(i_{L_1}(\ga_2)+\nu_{L_1}(\ga_2)\right)|\le
n-1.\lb{8.34} \eea By (\ref{8.32})-(\ref{8.34}), we have
 \bea
\left|\left(i_{L_0}(\ga^*)+\nu_{L_0}(\ga^*)\right)-
 \left(i_{L_1}(\ga^*)+\nu_{L_1}(\ga^*)\right)\right|
\le n-1,\nn\eea
 together with (\ref{8.31}), it implies Lemma 7.1.
\hfill\hb

Note that we can also prove Lemma 7.1 by Lemma 1.1, Proposition 6.1
and computation of the H${\rm \ddot{o}}$rmander index similarly as
the proof of Theorem 3.3 of \cite{LZZ}.

\noindent{\bf Lemma 7.2.} {\it Let $\ga\in \P_\tau(2n)$ be extended
to $[0,+\infty)$ by $\ga(\tau+t)=\ga(t)\ga(\tau)$ for all $t>0$.
Suppose $\ga(\tau)=M=P^{-1}(I_2\diamond \td{M})P$ with $\td{M}\in
\Sp(2n-2)$ and $i(\ga)\ge n$. Then we have
      \be i(\ga,2)+2S_{M^2}^+(1)-\nu(\ga,2)\ge n+2.\lb{8.35}\ee}
{\bf Proof.}
    The proof is similar to that of Lemma 4.1 in \cite{LLZ} (also
 Lemma 15.6.3 of \cite{Long1}). We write it down briefly.
 By (19) and (20) of the proof of Lemma 3 on p.349-350 in \cite{Long1}. We have
  \bea && i(\ga,2)+2S_{M^2}^+(1)-\nu(\ga,2)\nn\\
       &=& 2i(\ga)+2S_M^+(1)+\sum_{\theta\in
       (0,\pi)}(S_M^+(e^{\sqrt{-1}\theta})\nn\\
       &&-(\sum_{\theta\in
       (0,\pi)}(S_M^-(e^{\sqrt{-1}\theta})+(\nu(M)-S_M^-(1))+(\nu_{-1}(M)-S_M^-(-1)))\nn\\
       &\ge& 2n+2S_M^+(1)-n\nn\\
       &=&n+2S_M^+(1)\nn\\
        &\ge& n+2,\lb{8.36}\eea
  where in the last inequality we have used $\ga(\tau)=M=P^{-1}(I_2\diamond
  \tilde{M})P$ and the fact $S_{I_2}^+(1)=1$.

  \hfill\hb

\noindent{\bf Lemma 7.3.} {\it For any $(\tau,x)\in
\mathcal{J}_b(\Sg,2)$ and $m\in \N$, we have
\bea i_{L_0}(x,m+1)-i_{L_0}(x,m)&\ge& 1,\lb{8.37}\\
    i_{L_0}(x,m+1)+\nu_{L_0}(x,m+1)-1&\ge&
   i_{L_0}(x,m+1)>i_{L_0}(x,m)+\nu_{L_0}(x,m)-1.\lb{8.38}\eea}
{\bf Proof.} Let $\ga$ be the associated symplectic path of
$(\tau,x)$ and we extend $\ga$ to $[0,+\infty)$ by
$\gamma|_{[0,\frac{k\tau}{2}]}=\gamma^k$ with $\gamma^k$ defined in
(\ref{uvw}) for any $k\in\N$. By (\ref{8.5}) and (\ref{8.9}), for
any $m\in \N$ we have
       \be \nu_{L_0}(x,m)\ge 1, \qquad \forall m\in \N.\lb{8.39}\ee
Since $H_\Sg$ is strictly convex, $H_\Sg''(x(t))$ is positive for
all $t\in \R$. So by Theorem 5.1 and Lemma 5.1 of \cite{Liu2}(see
Theorem 2.4 in Section 2), we have
      \bea i_{L_0}(x,m+1)&=&
      \sum_{0< t<\frac{(m+1)\tau}{2}}\nu_{L_0}(\ga(t))\nn\\&\ge&
       \sum_{0<
       t\le\frac{m\tau}{2}}\nu_{L_0}(\ga(t))\nn\\&=&\sum_{0<
       t<\frac{m\tau}{2}}\nu_{L_0}(\ga(t))+\nu_{L_0}(\ga(\frac{m\tau}{2}))\nn\\&=&
       i_{L_0}(x,m)+\nu_{L_0}(x,m)\nn\\&>&i_{L_0}(x,m)+\nu_{L_0}(x,m)-1.\lb{8.40}\eea
 Thus  we get (\ref{8.37}) and (\ref{8.38})  from (\ref{8.39})
 and (\ref{8.40}). This proves Lemma 7.3.
 \hfill\hb

$\,$

 \noindent{\bf Proof of Theorem 1.1.}
 It is suffices to consider the case
 $^\#\tilde{\mathcal{J}}_b(\Sg)<+\infty$. Since $-\Sg=\Sg$, for
 $(\tau,x) \in \mathcal{J}_b(\Sg,2)$ we have
          \bea &&H_\Sg(x)=H_\Sg(-x),\lb{8.41}\\
               &&H_\Sg'(x)=- H_\Sg'(-x),\lb{8.42}\\
                &&H_\Sg''(x)= H_\Sg''(-x).\lb{8.43}\eea
So $(\tau,-x)\in \mathcal{J}_b(\Sg,2)$. By (\ref{8.43}) and the
definition of $\ga_x$ we have that
 \be \ga_x=\ga_{-x}.\lb{8.44}\ee
So we have
       \bea &&(i_{L_0}(x,m),\nu_{L_0}(x,m))=(i_{L_0}(-x,m),\nu_{L_0}(-x,m)),\nn\\
       &&(i_{L_1}(x,m),\nu_{L_1}(x,m))=(i_{L_1}(-x,m),\nu_{L_1}(-x,m)),\quad \forall m\in
       \N.\lb{8.45}\eea
 So we can write
  \be \td{\mathcal {J}}_b(\Sg,2)=\{[(\tau_j,x_j)]|
j=1,\cdots,p\}\cup\{[(\tau_k,x_k)],[(\tau_k,-x_k)]|k=p+1,\cdots,p+q\}.\lb{8.46}\ee
with $x_j(\R)=-x_j(\R)$ for $j=1,\cdots,p$ and $x_k(\R)\neq
-x_k(\R)$ for $k=p+1,\cdots,p+q$. Here we remind that $(\tau_j,x_j)$
has minimal period $\tau_j$ for $j=1,\cdots,p+q$ and
$x_j(\frac{\tau_j}{2}+t)=-x_j(t), \;t\in\R$ for $j=1,\cdots,p$.

 By Lemma 6.3 we have an integer $K\ge 0$ and an injection map
 $\phi: \N+K\to \mathcal
{V}_{\infty,b}(\Sg,2)\times \N$. By (\ref{8.45}), $(\tau_k,x_k)$ and
$(\tau_k,-x_k)$ have the same $(i_{L_0},\nu_{L_0})$-indices.
 So by Lemma 6.3,
 without loss of generality, we can further require that
       \be {\rm Im} (\phi)\subseteq \{[(\tau_k,x_k)]|k=1,2,\cdots,p+q\}\times
       \N.\lb{8.47}\ee
By the strict convexity of $H_\Sg$ and (\ref{6.11}), we have
        \be \hat{i}_{L_0}(x_k)>0,\quad k=1,2,\cdots,p+q.\lb{8.48}\ee
Applying  Theorem 1.5 and Remark 5.1 to  the following associated
symplectic paths
$$\ga_1,\;\cdots,\;\ga_{p+q},\; \ga_{p+q+1},\;\cdots,\;\ga_{p+2q}$$
of
$(\tau_1,x_1),\;\cdots,\;(\tau_{p+q},x_{p+q}),\;(2\tau_{p+1},x_{p+1}^2),\;\cdots,\;
     (2\tau_{p+q},x_{p+q}^2)$ respectively,
there exists a vector $(R,m_1,\cdots,m_{p+2q})\in \N^{p+2q+1}$ such
that $R>K+n$ and
   \bea &&i_{L_0}(x_k, 2m_k+1)=R+i_{L_0}(x_k),\lb{8.49}\\
        && i_{L_0}(x_k,2m_k-1)+\nu_{L_0}(x_k,2m_k-1)\nn\\
        &=&R-(i_{L_1}(x_k)+n+S_{M_k}^+(1)-\nu_{L_0}(x_k)),\lb{8.50}\eea
    for $k=1,\cdots,p+q,$ $M_k=\ga_k(\tau_k)$, and
     \bea &&i_{L_0}(x_k, 4m_k+2)=R+i_{L_0}(x_k,2),\lb{8.51}\\
         &&i_{L_0}(x_k,4m_k-2)+\nu_{L_0}(x_k,4m_k-2)\nn\\
         &=&R-(i_{L_1}(x_k,2)+n+S_{M_k}^+(1)-\nu_{L_0}(x_k,2)),\lb{8.52}\eea
      for $k=p+q+1,\cdots,p+2q$ and $M_k=\ga_k(2\tau_k)=\ga_k(\tau_k)^2$.

By Proposition 5.1 and the proof of Theorem 1.5, we also have \bea
i(x_k,
2m_k+1)&=&2R+i(x_k),\lb{8.53}\\
         i(x_k,2m_k-1)+\nu(x_k,2m_k-1)
        &=&2R-(i(x_k)+2S_{M_k}^+(1)-\nu(x_k)),\lb{8.54}\eea
    for $k=1,\cdots,p+q,$ $M_k=\ga_k(\tau_k)$, and
     \bea i(x_k, 4m_k+2)&=&2R+i(x_k,2),\lb{8.55}\\
        i(x_k,4m_k-2)+\nu(x_k,4m_k-2)
         &=&2R-(i(x_k,2)+2S_{M_k}^+(1)-\nu(x_k,2)),\lb{8.56}\eea
      for $k=p+q+1,\cdots,p+2q$ and $M_k=\ga_k(2\tau_k)$.

From (\ref{8.47}), we can set
 \bea \phi(R-(s-1))=([(\tau_{k(s)}, x_{k(s)})],m(s)),\qquad
    \forall s\in S:=\left\{1,2,\cdots,\left[\frac{n}{2}\right]+1\right\},\lb{8.57}\eea
where $k(s)\in \{1,2,\cdots,p+q\}$ and $m(s)\in \N$.

We continue our proof  to study the symmetric and asymmetric orbits
separately. Let \be S_1=\{s\in S|k(s)\le p\},\qquad S_2=S\setminus
S_1.\lb{8.58}\ee
 We shall prove that
$^\#S_1\le p$ and $^\#S_2\le 2q$, together with the definitions of
$S_1$ and $S_2$, these yield Theorem 1.1.

\noindent{\it Claim 1.} $^\#S_1\le p$.

\noindent {\it Proof of Claim 1.} By the definition of $S_1$,
$([(\tau_{k(s)},
 x_{k(s)})],m(s))$ is symmetric when $k(s)\le p$. We further prove
 that $m(s)=2m_{k(s)}$ for $s\in S_1$.

  In fact, by the definition of $\phi$ and Lemma 6.3, for all $s=1,2,\cdots,\left[\frac{n}{2}\right]+1$ we have
         \bea  i_{L_0}(x_{k(s)},m(s))&\le & (R-(s-1))-1=R-s \nn\\
         &\le &
         i_{L_0}(x_{k(s)},m(s))+\nu_{L_0}(x_{k(s)},m(s))-1.\lb{8.59}\eea
 By the strict convexity of $H_\Sg$, from Theorem 2.4, we have $i_{L_0}(x_{k(s)})\ge 0$, so there holds
   \bea i_{L_0}(x_{k(s)},m(s))\le R-s< R\le R+i_{L_0}(x_{k(s)})=i_{L_0}(x_{k(s)},2m_{k(s)}+1),\lb{8.60}\eea
 for every $s=1,2,\cdots,\left[\frac{n}{2}\right]+1$, where we have used
 (\ref{8.49}) in the last equality. Note that the proofs of (\ref{8.59}) and
 (\ref{8.60}) do not depend on the condition $s\in S_1$.

By Lemma 1.2, we have
   \be
i_{L_1}(x_k)+S_{M_k}^+(1)-\nu_{L_0}(x_k)\ge \frac{1-n}{2},\quad
\forall k=1,\cdots,p.\lb{8.64}\ee
 Also for $1\le
s\le \left[\frac{n}{2}\right]+1$, we have
   \be -\frac{n+3}{2}<-(1+\frac{n}{2})\le -(\left[\frac{n}{2}\right]+1)\le
   -s.\lb{8.65}\ee
Hence by (\ref{8.59}),(\ref{8.64}) and(\ref{8.65}), if $k(s)\le p$
we have
 \bea
&&i_{L_0}(x_{k(s)},2m_{k(s)}-1)+\nu_{L_0}(x_{k(s)},2m_{k(s)}-1)-1\nn\\
&=&
R-(i_{L_1}(x_{k(s)})+n+S_{M_{k(s)}}^+(1)-\nu_{L_0}(x_{k(s)}))-1\nn\\
&\le&R-\frac{1-n}{2}-1-n=R-\frac{n+3}{2}<R-s\nn\\
&\le&
i_{L_0}(x_{k(s)},m(s))+\nu_{L_0}(x_{k(s)},m(s))-1.\lb{8.66}\eea
 Thus
by (\ref{8.60}) and (\ref{8.66}) and Lemma 7.3 we have
 \be 2m_{k(s)}-1< m(s)<2m_{k(s)}+1.\lb{8.67}\ee
 Hence
 \be m(s)=2m_{k(s)}.\lb{8.68}\ee
So we have
 \be \phi(R-s+1)=([(\tau_{k(s)},x_{k(s)})],2m_{k(s)}),\qquad \forall
 s\in S_1.\lb{8.69}\ee
Then by the injectivity of $\phi$, it induces another injection map
 \be \phi_1:S_1\rightarrow \{1,\cdots,p\}, \;s\mapsto k(s).\lb{8.70}\ee
 There for $^\#S_1\le p$. Claim 1 is proved.

 \noindent{\it Claim 2.} $^\#S_2\le 2q$.

\noindent{\it Proof of Claim 2.} By the  formulas
(\ref{8.53})-(\ref{8.56}), and (59) of \cite{LLZ} (also Claim 4 on
p. 352 of \cite{Long1}), we have \be m_k=2m_{k+q}\quad {\rm for}\;\;
k=p+1,p+2,\cdots,p+q.\lb{8.71}\ee
 We set $\mathcal {A}_k=i_{L_1}(x_k,2)+S_{M_k}^+(1)-\nu_{L_0}(x_k,2)$
and $\mathcal {B}_k=i_{L_0}(x_k,2)+S_{M_k}^+(1)-\nu_{L_1}(x_k,2)$,
$p+1\le k\le p+q$, where $M_k=\ga_k(2\tau_k)=\ga(\tau_k)^2$. By
(\ref{6.10}),  we have
 \be
\mathcal {A}_k+\mathcal
{B}_k=i(x_k,2)+2S_{M_k}^+(1)-\nu(x_k,2)-n,\;\;\;p+1\le k\le p+q
.\lb{8.72}\ee By similar discussion of the proof of Lemma 1.1, for
any $p+1\le k\le p+q$ there exist $P_k\in \Sp(2n)$ and $\td{M}_k\in
\Sp(2n-2)$ such that \be \ga(\tau_k)=P_k^{-1}(I_2\diamond
\td{M}_k)P_k.\lb{8.73}\ee
 Hence by Lemma 7.2 and (\ref{8.72}), we
have \be \mathcal {A}_k+\mathcal {B}_k\ge n+2-n=2.\lb{8.74}\ee By
Theorem 2.3, there holds
 \bea
|\mathcal {A}_k-\mathcal
{B}_k|&=&|(i_{L_0}(x_k,2)+\nu_{L_0}(x_k,2))-(i_{L_1}(x_k,2)+\nu_{L_1}(x_k,2))|\le
n.\lb{8.75}\eea
 So by (\ref{8.74}) and (\ref{8.75}) we have \be
\mathcal {A}_k\ge \frac{1}{2}((\mathcal {A}_k+\mathcal
{B}_k)-|\mathcal {A}_k-\mathcal {B}_k|)\ge \frac{2-n}{2},\quad
p+1\le k\le p+q.\lb{8.76}\ee
 By
(\ref{8.52}), (\ref{8.59}), (\ref{8.65}), (\ref{8.71}) and
(\ref{8.76}), for $p+1\le k(s)\le p+q$ we have
\bea &&i_{L_0}(x_{k(s)},2m_{k(s)}-2)+\nu_{L_0}(x_{k(s)},2m_{k(s)}-2)-1\nn\\
     &=&i_{L_0}(x_{k(s)},4m_{k(s)+q}-2)+\nu_{L_0}(x_{k(s)},4m_{k(s)+q}-2)-1\nn\\
     &=&R-(i_{L_1}(x_{k(s)},2)+n+S_{M_{k(s)}}^+(1)-\nu_{L_0}(x_{k(s)},2))-1\nn\\
     &=&R-\mathcal{A}_{k(s)}-1-n\nn\\
     &\le&R- \frac{2-n}{2}-1-n\nn\\
     &=& R-(2+\frac{n}{2})\nn\\
     &<& R-s\nn\\
     &\le&
     i_{L_0}(x_{k(s)},m(s))+\nu_{L_0}(x_{k(s)},m(s))-1.\lb{8.77}\eea
Thus by (\ref{8.60}), (\ref{8.77}) and Lemma 7.3, we have
     \be 2m_{k(s)}-2<m(s)<2m_{k(s)}+1,\qquad p<k(s)\le p+q.\lb{8.78}\ee
So \be m(s)\in \{2m_{k(s)}-1,2m_{k(s)}\}, \qquad {\rm
for}\;\;p<k(s)\le p+q.\lb{8.79}\ee
 Especially this yields that for any $s_0$ and $s\in
S_2$, if $k(s)=k(s_0)$, then
 \be m(s)\in
\{2m_{k(s)}-1,2m_{k(s)}\}=\{2m_{k(s_0)}-1,2m_{k(s_0)}\}.\lb{8.80}\ee
Thus by the injectivity of the map $\phi$ from Lemma 3.3, we have
\be ^\#\{s\in S_2|k(s)=k(s_0)\}\le 2.\lb{8.81}\ee This yields Claim
2.

     By Claim 1 and Claim 2, we have
     \bea
     ^\#\td{\mathcal{J}}_b(\Sg)=^\#\td{\mathcal{J}}_b(\Sg,2)=p+2q\ge
     ^\#S_1+^\#S_2
     =\left[\frac{n}{2}\right]+1.\lb{8.82}\eea
The proof of Theorem 1.1 is complete. \hfill\hb
\setcounter{equation}{0}
\section {Proof of Theorem 1.2.}

{\bf Proof of Theorem 1.2.} We prove Theorem 1.2 in three steps.

 \noindent{\it Step 1.} Applying  Theorem 1.5.

 If
$^\#\tilde{\mathcal{J}}_b(\Sg)<+\infty$, we write
  \bea \td{\mathcal {J}}_b(\Sg,2)=\{[(\tau_j,x_j)]|
j=1,\cdots,p\}\cup\{[(\tau_k,x_k)],[(\tau_k,-x_k)]|k=p+1,\cdots,p+q\},
\nn\eea where $(\tau_j,x_j)$ is symmetric with minimal period
$\tau_j$ for $j=1,\cdots,p$, and $(\tau_k,x_k)$ is asymmetric with
minimal period $\tau_k$ for $k=p+1,\cdots,p+q$, for simplicity we
have set $q=\mathfrak{A}(\Sg)$ with $\mathfrak{A}(\Sg)$ defined in
Theorem 1.2.

 By Lemma 6.3, there exist $0\le K\in \Z$ and injection map $\phi: \N+K\to \mathcal
{V}_{\infty,b}(\Sg,2)\times \N$ such that (i) and (ii) in Lemma 6.3
hold. By the same reason for (\ref{8.47}), we can require that
       \be {\rm Im} (\phi)\subseteq \{[\tau_k,x_k)]|k=1,2,\cdots,p+q\}\times
       \N.\lb{8.47'}\ee
Set $r=p+q$.
 By (\ref{8.48}) we have ${\hat{i}}_{L_0}(x_j)>0$ for $j=1,\cdots,r$.  Applying Theorem
 1.5 and Remark 5.1 to the collection of symplectic paths $\ga_1,\,\ga_2,\,\cdots,\,\ga_r$,
 there exists a vector $(R, m_1, m_2,\cdots,m_r)\in \N^{r+1}$
 such that $R>K+n$ and
   \bea &&\nu_{L_0}(\ga_j, 2m_j\pm 1)=\nu_{L_0}(\ga_k),\lb{9.1}\\
   && i_{L_0}(\ga_j, 2m_j-1)+\nu_{L_0}(\ga_j,2m_k-1)=R-(i_{L_1}(\ga_j)+n+S_{M_j}^+(1)-\nu_{L_0}(\ga_j)),\lb{9.2}\\
&&i_{L_0}(\ga_j,2m_k+1)=R+i_{L_0}(\ga_j),\lb{9.3}\eea where $\ga_j$
is the associated symplectic path of
  $(\tau_j,x_j)$ and $M_j=\ga_j(\tau_j)$, $1\le j\le r$.

\noindent{\it Step 2.} We prove that \bea K_1:=\min
  \{i_{L_1}(\ga_j)+S_{M_j}^+(1)-\nu_{L_0}(\ga_j)|j=1,\cdots,r\}\ge 0.\lb{9.4}\eea
  By the strict convexity of $H_\Sg$, Theorem 2.4 yields
  \be i_{L_1}(\ga_j)\ge 0.\lb{9.5}\ee
  By the
  nondegenerate assumption in Theorem 1.2 we have $\nu_{L_0}(\ga_j, m)=1$ for
  $ 1\le j\le r,\;m\in \N$. By  similar discussion of Lemma 1.1,
  there exist $P_j\in \Sp(2n)$ and $\tilde {M}_j\in \Sp(2n-2)$ such
  that
        $$M_j=P_j^{-1}(I_2\diamond \tilde{M}_j)P_j.$$
 So we have
       \be S_{M_j}^+(1)=S_{I_2\diamond \tilde{M}_j}^+(1)=S_{I_2}^+(1)
       +S_{\td{M}_j}^+(1)\ge S_{I_2}^+(1)=1.\lb{9.6}\ee
Thus (\ref{9.5}) and (\ref{9.6}) yield
       $$K_1\ge 0.$$

\noindent{\it Step 3.} Complete the proof of Theorem 1.2.

By (\ref{8.47'}), we set
$\phi(R-(s-1))=([(\tau_j(s),x_{j(s)})],m(s))$ with $j(s)\in
\{1,\cdots,r\}$ and $m(s)\in \N$ for $s=1,\cdots,n$. By Lemma 6.2 we
have
  $$i_{L_0}(x_{j(s)},m(s))\le R-(s-1)-1=R-s\le i_{L_0}(x_{j(s)},m(s))+\nu_{L_0}(x_{j(s)},m(s))-1.$$
By (\ref{9.2}) and (\ref{9.4}) for $s=1,\cdots,n$,
     \bea &&i_{L_0}(x_{j(s)},2m_{j(s)}-1)+\nu_{L_0}(x_{j(s)},2m_{j(s)}-1)-1\le
     R-K_1-1-n<R-n\nn\\
     &&\le R-s\le i_{L_0}(x_{j(s)},m(s))+\nu_{L_0}(x_{j(s)},m(s))-1.\nn\eea
By (\ref{8.38}), we have
      $$2m_{j(s)}-1<m(s), \quad s=1,\cdots,n.$$
For $s=1,\cdots,n$, there holds
      $$i_{L_0}(x_{j(s)}, m(s))\le R-s<R\le
      i_{L_0}(x_{j(s)},2m_{j(s)}+1),$$
then by (\ref{8.38}), we have
      $$m(s)<2m_{j(s)}+1, \quad s=1,\cdots,n.$$
 Thus
          \bea m(s)=2m_{j(s)},\quad \quad s=1,\cdots,n.\lb{9.7}\eea
By (ii) of Lemma 6.3 again, if $s_1\neq s_2$, we have $m(s_1)\neq
m(s_2)$. By (\ref{9.7}) we have $j(s_1)\neq j(s_2)$. So $j(s)'s$ are
mutually different for $s=1,\cdots,n$. Since $j(s)\in
\{1,2,\cdots,r\}$, we have
     $$r\ge n.$$
Hence \be
^\#\tilde{\mathcal{J}}_b(\Sg)=^\#\tilde{\mathcal{J}}_b(\Sg,2)=
p+2q=r+q\ge n+q=n+\mathfrak{A}(\Sg).\lb{9.8}\ee The proof of Theorem
1.2 is complete. \hfill\hb

\noindent{\it Acknowledgments.} The authors thank Professor Y. Long
for stimulating and very useful discussions, and for encouraging us
to study Maslov-type index theory and its iteration theory. They
also thank Professor C. Zhu for his valuable suggestions.

\bibliographystyle{abbrv}

\end{document}